\documentclass{gtart_h}

\def\ifplaintex{\expandafter\ifx\csname documentclass\endcsname\relax}

\def\gtp{{\mathsurround=0pt\it $\cal G\mskip-2mu$eometry \&\ 
$\cal T\!\!$opology $\cal P\!$ublications}}  

\def\recd{{\small Received:\qua\receiveddate\ifx\reviseddate\relax
\else\qquad Revised:\qua\reviseddate\fi\par}} 

\def\lognumber#1{\def\thelognumber{#1}}
\def\volumenumber#1{\def\thevolumenumber{#1}}
\def\volumeyear#1{\def\thevolumeyear{#1}}
\def\papernumber#1{\def\thepapernumber{#1}}
\def\pagenumbers#1#2{\def\startpage{#1}\def\finishpage{#2}}
\def\published#1{\def\publishdate{#1}}

\def\received#1{\def\receiveddate{#1}}
\def\revised#1{\def\reviseddate{#1}}
\def\accepted#1{\def\accepteddate{#1}}

\def\asciiauthors#1{\def\theasciiauthors{#1}}
\def\asciiaddress#1{\def\theasciiaddress{#1}}
\def\asciiemail#1{\def\theasciiemail{#1}}

\def\coverauthors#1{\def\thecoverauthors{#1}}
\long\def\asciiabstract#1{\long\def\theasciiabstract{#1}}

\let\\\par\let\thelognumber\relax\let\thevolumenumber\relax
\let\thepapernumber\relax\let\thevolumeyear\relax\let\startpage\relax
\let\finishpage\relax\let\publishdate\relax\let\receiveddate\relax
\let\reviseddate\relax\let\accepteddate\relax\let\theasciititle\relax
\let\theasciiauthors\relax\let\theasciiaddress\relax
\let\theasciiabstract\relax

\let\thecoverauthors\relax\let\theasciiemail\relax

\ifplaintex
\font\logobig=cmssbx10 scaled 3836
\font\logomed=cmssbx10 scaled 2557
\else
\font\logobig=cmssbx10 scaled 4200
\font\logomed=cmssbx10 scaled 2800
\fi

\long\def\makeagttitle{   
\count0=\startpage
\agt\hfill      
\hbox to 45truept{\vbox to 0pt{\vglue -13truept{\logomed A\kern -.37em{\logobig 
T}\kern -.38em G}\vss}\hss}
\break
{\small Volume \thevolumenumber\ (\thevolumeyear)
\startpage--\finishpage\nl
Published: \publishdate}

\vglue .25truein

{\parskip=0pt\leftskip 0pt plus
1fil\def\\{\par\smallskip}{\Large\bf\thetitle}\par\medskip} \vglue
0.05truein

{\parskip=0pt\leftskip 0pt plus 1fil\def\\{\par}{\sc\theauthors}
\par\medskip}%
 
\vglue 0.03truein 

{\small\leftskip 25truept\rightskip 25truept{\bf Abstract}\stdspace\theabstract

{\bf AMS Classification}\stdspace\theprimaryclass
\ifx\thesecondaryclass\relax\else; \thesecondaryclass\fi\par
{\bf Keywords}\stdspace \thekeywords\par}\vglue 7truept

}   

\ifplaintex
\font\phead=cmsl9 scaled 950
\font\pnum=cmbx10 scaled 913
\font\pfoot=cmsl9 scaled 950
\headline{\vbox to 0pt{\vskip -4.5mm\line{\small\phead\ifnum
\count0=\startpage ISSN 1472-2739 (on-line) 1472-2747 (printed)
\hfill {\pnum\folio}\else\ifodd\count0\def\\{ }%
\ifx\theshorttitle\relax\thetitle\else\theshorttitle\fi\hfill{\pnum\folio}
\else\def\\{ and }{\pnum\folio}\hfill\ifx\theshortauthors\relax\theauthors
\else\theshortauthors\fi\fi\fi}\vss}}
\footline{\vbox to 0pt{\vglue 0mm\line{\small\pfoot\ifnum\count0=\startpage
\copyright\ \gtp\hfill\else
\agt, Volume \thevolumenumber\ (\thevolumeyear)\hfill\fi}\vss}}
\else
\font=cmsl9 scaled 1050
\font\lnum=cmbx10 
\font=cmsl9 scaled 1050
\makeatletter
\def\@oddhead{{\small\lhead\ifnum\count0=\startpage ISSN 1472-2739 
(on-line) 1472-2747 (printed)\hfill {\lnum\number\count0}\else\ifodd\count0
\def\\{ }\ifx\theshorttitle\relax \thetitle \else\theshorttitle\fi\hfill
{\lnum\number\count0}\else\def\\{ and }{\lnum\number\count0}
\hfill\ifx\theshortauthors\relax 
\theauthors\else\theshortauthors\fi\fi\fi}}\def\@evenhead{\@oddhead}
\def\@oddfoot{\small\lfoot\ifnum\count0=\startpage\copyright\ \gtp\hfill\else
\agt, Volume \thevolumenumber\ (\thevolumeyear)\hfill\fi}
\def\@evenfoot{\@oddfoot}
\makeatother
\fi
\let\maketitlepage\makeagttitle
\let\makeshorttitle\maketitlepage
\let\maketitle\maketitlepage

\newwrite\gtoutfile
\long\gdef\makeheadfile{  
{\def\\{, }\def\s{ }
\immediate\openout\gtoutfile head.xxx
\immediate\write\gtoutfile{Proxy-for: \ifx\theasciiauthors\relax
\theauthors\else\theasciiauthors\fi\s<\ifx\theasciiemail\relax\theemail\else\theasciiemail\fi>}
\immediate\write\gtoutfile{\noexpand\\}
\immediate\write\gtoutfile{Authors: \ifx\theasciiauthors\relax
\theauthors\else\theasciiauthors\fi}
{\def\\{ }\immediate\write\gtoutfile{Title: \ifx\theasciititle\relax
\thetitle\else\theasciititle\fi}}
\immediate\write\gtoutfile{Subj-class: GT or SG, GR etc}
\immediate\write\gtoutfile{MSC-class: \theprimaryclass\ifx\thesecondaryclass\relax\else, \thesecondaryclass\fi}
\immediate\write\gtoutfile{Journal-ref: Algebr. Geom. Topol. \thevolumenumber\s
(\thevolumeyear) \startpage-\finishpage}
\immediate\write\gtoutfile{Comments: Published by Algebraic and
Geometric Topology at}
\immediate\write\gtoutfile{\s\s\s  http://www.maths.warwick.ac.uk/agt/AGTVol\thevolumenumber/agt-\thevolumenumber-\thepapernumber.abs.html}
\immediate\write\gtoutfile{\noexpand\\}
\immediate\write\gtoutfile{}
\ifx\theasciiabstract\relax
\immediate\write\gtoutfile{\theabstract}\else
\immediate\write\gtoutfile{\theasciiabstract}\fi
\immediate\write\gtoutfile{}
\immediate\write\gtoutfile{\noexpand\\}
\immediate\write\gtoutfile{}
\immediate\closeout\gtoutfile}}  

\def\maketitlepage{\makeagttitle\makeheadfile}
\let\makeshorttitle\maketitlepage
\let\maketitle\maketitlepage

\lognumber{41}
\volumenumber{5}
\volumeyear{2005}
\papernumber{41}
\pagenumbers{999}{1026}
\received{1 March 2005} 
\revised{28 June 2005}
\accepted{18 July 2005}
\published{18 August 2005}

\usepackage{amsmath}
\usepackage{amssymb}
\usepackage{graphics}

\let\tilde\widetilde
\def\figref#1{\hyperlink{#1anchor}{Figure~\ref*{#1}}}
\def\tableref#1{\hyperlink{#1anchor}{Table~\ref*{#1}}}
\def\anchor#1{\noindent\hypertarget{#1anchor}{\smash{$\phantom{99}$}}}
\allowdisplaybreaks

\newcommand\hiv{\mathbf H^4}
\newcommand\hn{\mathbf H^n}

\newcommand\rnm{\mathbf R^{n-1}}

\newcommand\riv{\mathbf R^4}
\newcommand\z{\mathbf Z}

\newcommand\mbar{\overline M}
\newcommand\ncl{\left<\left<}
\newcommand\ncr{\right>\right>}

\DeclareMathOperator{\Isom}{Isom}

\DeclareMathOperator{\stab}{stab}

\newcommand{\arrtop}[2]{\xrightarrow[#2]{#1}}

\newcommand\bd{\partial}

\begin{document}

\newtheorem{theorem}{Theorem}[section]
\newtheorem{proposition}[theorem]{Proposition}
\theoremstyle{definition}
\newtheorem{remark}[theorem]{Remark}
\newtheorem{example}[theorem]{Example}

\title{Complements of tori and Klein bottles in the\\4-sphere that have hyperbolic structure}
\shorttitle{Complements of tori and Klein bottles with hyperbolic structure}

\author{Dubravko Ivan\v si\'c\\John G. Ratcliffe\\Steven T. Tschantz}
\asciiauthors{Dubravko Ivansic, John G. Ratcliffe and Steven T. Tschantz}
\shortauthors{Ivan\v si\'c, Ratcliffe and Tschantz}
\coverauthors{Dubravko Ivan\noexpand\v si\noexpand\'c\\John G. Ratcliffe\\Steven T. Tschantz}
\address{{\rm DI:}\qua Department of Mathematics and Statistics, Murray State
University\\Murray, KY 42071, USA\\\smallskip\\{\rm JR and ST:}\qua Department of Mathematics, Vanderbilt University\\Nashville, TN 37240, USA}

\asciiaddress{DI: Department of Mathematics and Statistics\\Murray State
University, Murray, KY 42071, USA\\JR and ST: Department of Mathematics, Vanderbilt University\\Nashville, TN 37240, USA}

\asciiemail{dubravko.ivansic@murraystate.edu, ratclifj@math.vanderbilt.edu,
tschantz@math.vanderbilt.edu}
\gtemail{\mailto{dubravko.ivansic@murraystate.edu}, 
\mailto{ratclifj@math.vanderbilt.edu},
\mailto{tschantz@math.vanderbilt.edu}}

\begin{abstract}
Many noncompact hyperbolic 3-manifolds are topologically\break complements
of links in the 3-sphere.  Generalizing to dimension 4, we construct a
dozen examples of noncompact hyperbolic 4-manifolds, all of which are
topologically complements of varying numbers of tori and Klein bottles
in the 4-sphere.  Finite covers of some of those manifolds are then
shown to be complements of tori and Klein bottles in other
simply-connected closed 4-manifolds.  All the examples are based on a
construction of Ratcliffe and Tschantz, who produced 1171 noncompact
hyperbolic 4-manifolds of minimal volume.  Our examples are finite
covers of some of those manifolds.
\end{abstract}

\asciiabstract{%
Many noncompact hyperbolic 3-manifolds are topologically complements
of links in the 3-sphere.  Generalizing to dimension 4, we construct a
dozen examples of noncompact hyperbolic 4-manifolds, all of which are
topologically complements of varying numbers of tori and Klein bottles
in the 4-sphere.  Finite covers of some of those manifolds are then
shown to be complements of tori and Klein bottles in other
simply-connected closed 4-manifolds.  All the examples are based on a
construction of Ratcliffe and Tschantz, who produced 1171 noncompact
hyperbolic 4-manifolds of minimal volume.  Our examples are finite
covers of some of those manifolds.}


\maketitle

\section{Introduction}
\label{introduction}

Let $\hn$ be the $n$-dimensional hyperbolic space and let  $G$ be a discrete
subgroup of $\Isom\hn$, the isometries of $\hn$.
If $G$ is torsion-free, then $M=\hn/G$ is a hyperbolic manifold of
dimension $n$.
A hyperbolic manifold of interest in this paper is also complete, noncompact
and has finite volume and the term ``hyperbolic manifold'' will be understood to include
those additional properties.
Such a manifold $M$ is the interior of a compact manifold with boundary
$\mbar$.
Every boundary component $E$ of $\mbar$ is a compact flat (Euclidean) manifold,
i.e. a manifold of the form $\rnm/K$, where $K$ is a discrete subgroup of
$\Isom\rnm$, the isometries of $\rnm$.
For simplicity, we sometimes inaccurately call $E$ a boundary component of $M$ (rather than
$\mbar$, whose boundary component it really is).

When $n=3$, $M$ is the interior of a compact 3-manifold $\mbar$ whose boundary consists
of tori and Klein bottles.
When all the boundary components are tori, $M$ may be diffeomorphic to $S^3$ with several
closed solid tori removed, which in turn is diffeomorphic to $S^3-A$, where $A$ is the link
consisting of the core circles of the solid tori.
Indeed, it is a well known fact that many hyperbolic 3-manifolds are link complements in $S^3$.

The first example of a generalization of this situation to dimension 4 was given by
Ivan\v si\'c~\cite{Ivansic3}, where it was shown that $\tilde M_{1011}$, the orientable double
cover of one of the 1049 nonorientable hyperbolic 4-manifolds that Ratcliffe and Tschantz constructed 
\cite{Ratcliffe-Tschantz}, is a complement of 5 tori in $S^4$.
Finite covers of  $\tilde M_{1011}$ were also shown to be complements of a collection
of tori in a simply-connected 4-manifold with even Euler characteristic.

In this paper, we find 11 additional examples of hyperbolic 4-manifolds that are complements in
$S^4$ of a set  of tori or a combination of tori and Klein bottles.
The list of examples is in Theorem~\ref{mainthm}.
Like $\tilde M_{1011}$, all of them are orientable double covers of some of Ratcliffe and
Tschantz's nonorientable manifolds.
Proving that those manifolds are complements in $S^4$ is typically trickier than in \cite{Ivansic3},
since the many symmetries of $M_{1011}$ enabled an easier computation.

By taking finite covers of some of the link complements we also find additional examples of
hyperbolic 4-manifolds that are complements of a collection of tori and Klein bottles in a
simply-connected 4-manifold with even Euler characteristic.
Example~\ref{m1091supp} displays 8 different families of examples, illustrating the richness
of such objects, many more of which likely exist.

Manifolds $\tilde M_{1091}$ and  $\tilde M_{1011}$ were used by
Ratcliffe and Tschantz~\cite{Ratcliffe-Tschantz2} to construct the first examples (infinitely many) of
aspherical homology 4-spheres.
Three other examples in this paper can be used for the same purpose, see Example~\ref{aspherical}. 

When looking for examples of link complements in $S^4$ among the orientable double covers
of the 1049 nonorientable manifolds, the main difficulty is to reduce the number of potential examples
with which to experiment.
We use a homological and a group-theoretical criterion (see \S\ref{examples}) to rule out all
except 49 manifolds from having the desired property.
Then we show that some of the remaining 49 manifolds have orientable double covers that are
complements of tori and Klein bottles in the 4-sphere.
Essentially, we use group presentations to show that the manifold $N$ resulting from closing up the
boundary components of the double covers is simply-connected.
This, along with the fact that the Euler characteristic $\chi(N)$ is~2, guarantees that $N$ is
homeomorphic to $S^4$.

Ratcliffe and Tschantz constructed their manifolds by side-pairings of a hyperbolic noncompact
right-angled 24-sided polyhedron, making computations with them laborious, especially if they
have to be done repeatedly.
The two criteria were checked using a computer, but we emphasize that the results of this
paper do not depend on a computer calculation, since it is verified by hand that some of the
``promising'' examples indeed have a double cover that is a complement in the 4-sphere.

\section{Some preliminary facts}
\label{prelims}

Let $M$ be hyperbolic $n$-manifold.
We say that $M$ {\it is a (codimension-2) complement in} $N$ if $M=N-A$, where $N$ is a closed
$n$-manifold and $A$ is a closed $(n-2)$-submanifold of $N$ that has a tubular neighborhood in
$N$ and has as many components as $\bd\mbar$.
Then every component $E$ of $\bd\mbar$, being the boundary of a tubular neighborhood, must be an
$S^1$-bundle over a component of $A$.
The $(n-1)$-manifold $E$ is flat.
The following theorem (see \cite{Ivansic2} or \cite{Apanasov} Theorem~6.40) summarizes when a
flat manifold is an $S^1$-bundle and, when it is, what its $S^1$-fibers are.

\begin{theorem}
\label{bundlesref}
Let $E=\rnm/K$ be a compact flat $(n-1)$-manifold, where $K$ is a discrete subgroup
of $\Isom\rnm$.
Then $E$ is an $S^1$-bundle over a manifold $B$ if and only if there exists a translation
$t\in K$ such that $\left<t\right>$ is a normal subgroup of $K$ and $t$ is not
a power of any element of $K$ other than $t^{\pm 1}$.
If the translation is given by $x\mapsto x+v$ and any element of~$K$ by $x\mapsto Ax+a$,
where $A\in O(n-1)$, $a\in\rnm$, the normality of $\left<t\right>$ can
be expressed as $Av=\pm v$ for every $A$ such that $Ax+a$ is an element of $K$.
Furthermore, if $n-1\ne 4,5$ (we are interested in $n-1=3$ in this paper) the manifold
$B$ above is a flat manifold.
\qed
\end{theorem}

We say that a translation $t$ is {\it normal in $E$} if $\left<t\right>$ is a normal subgroup, we say
it is {\it primitive in $E$} if it is not a power of any other element of $K$.

By Theorem~\ref{bundlesref}, if $M$ is a hyperbolic 4-manifold that is a complement in $N$,
it must be a complement of flat 2-manifolds, that is, tori and Klein bottles.
We are interested in the case when $N=S^4$.
When $M=N-A$, an immediate consequence of the fact that components of $A$ are flat is that
$\chi(M)=\chi(N)$, so examples of complements in $S^4$ must be searched for among manifolds
with $\chi(M)=2$.
Nonorientable examples of Ratcliffe and Tschantz are very suitable for this
purpose, since their Euler characteristic is 1: passing to the orientable double cover will
give us the required Euler characteristic.

Once we have found a hyperbolic manifold $M$ all of whose boundary components are
\makebox{$S^1$-bundles}, $M=N-A$.
To show that $N$ is a topological $S^4$, it suffices to show that $\pi_1 N=1$ and $\chi(N)=2$
(see \cite{Gompf-Stipsicz} or \cite{Freedman-Quinn}).
The following application of van Kampen's theorem shows how to compute $\pi_1 N$ from
$\pi_1 M$; details are in \cite{Ivansic3}.

\begin{proposition}
\label{presentationofN}
Let $M$ be a complement in $N$ and let $\bd\mbar=E_1\cup\dots\cup E_m$.
If $t_1,\dots,t_m\in \pi_1 M$ represent fibers of the $S^1$-bundles $E_1,\dots,E_m$
then $\pi_1 N = \pi_1 M/ \ncl t_1,\dots,t_m \ncr$, where $\ncl A \ncr$ denotes the normal
closure of a subset $A\subset \pi_1 M$.
In other words, if we have a presentation for $M$, the presentation for $N$ is
obtained by adding relations $t_1=1,\dots,t_m=1$.
\qed
\end{proposition}

Hyperbolic manifolds $M$ of interest in this paper are given by side-pairings of a hyperbolic
finite-volume noncompact polyhedron $Q$.
Side-pairings determine sets of points on the polyhedron, called {\it cycles}, that are identified in
the manifold.
Boundary components of $\mbar$ correspond to cycles of ideal vertices of $Q$.
If $[v_i]$ is the cycle of an ideal vertex corresponding to $E_i$, then $\stab v_i$ is
isomorphic to $\pi_1 E_i$, making the inclusion $\pi_1 E_i\to\pi_1 M$
injective.
Elements of $\pi_1 E_i$ are parabolic isometries, thus, they preserve horospheres centered at
$v_i$ and act like Euclidean isometries on them.

Let $\Gamma$ be a discrete subgroup of $\Isom\hn$ with fundamental polyhedron $P$.
Then $\Gamma$ is generated by side-pairings of $P$, isometries that send a side of $P$ to another
side of $P$.
If $s$ is an isometry sending side $S$ of $P$ to $S'$, then, if we start moving from $P$ and
pass through~$S$, we wind up in $s^{-1}P$.
Similarly, if we take a path $P$ to $\gamma P$ that passes through interiors of translates of sides
$S_1,\dots,S_m$ whose side-pairings are $s_1,\dots,s_k$, then $\gamma=s_1^{-1}\dots s_m^{-1}$.
Furthermore, if $G$ is a finite-index subgroup of $\Gamma$ with transversal $X$ (that is,
a set of right-coset representatives, so $\Gamma=\cup_{x\in X} Gx$), then the fundamental
polyhedron of $G$ is $Q=\cup_{x\in X} xP$, where elements of $X$ can be chosen so
that $Q$ is connected.

The fundamental groups $G$ of Ratcliffe and Tschantz's manifolds are all minimal index
torsion-free subgroups of a certain reflection group $\Gamma$, which is generated by reflections in a
10-sided right-angled polyhedron $P$.
The side-pairings generating $\Gamma$ are, of course, reflections in the sides of $P$.
Every subgroup $G$ corresponding to a Ratcliffe-Tschantz manifold has the same transversal in
$\Gamma$,  a finite group $K$ isomorphic to $\z_2^4$ and generated by reflections in the
four coordinate planes in the ball model of $\hiv$.
The fundamental polyhedron of $G$ is $Q=KP$, a 24-sided regular right-angled hyperbolic polyhedron
that we call the {\it 24-cell}.
Let $f$ be a side-pairing of $Q$ that sends side $R$ to side $R'$.
As was shown in \cite{Ratcliffe-Tschantz} and \cite{Ivansic3}, $f$ has to
have form $f=kr=r'k$, where $k\in K$ is an element such that $k(R)=R'$, and $r$, $r'$ are reflections
in $R$, $R'$, respectively.
We will call $k$ the {\it $K$-part of $f$}.
Thus, it is enough to specify elements $k\in K$ in order to define a side-pairing of $Q$.

As is clear from Theorem~\ref{bundlesref} and Proposition~\ref{presentationofN}, we will need
the ability to recognize when a composite of side-pairings is a translation.
The following proposition, proved in \cite{Ivansic3}, will be used for that purpose:
 
\begin{proposition}
\label{moving}
Let $\Gamma$ be generated by reflections in the sides of a polyhedron $P$ and let $G$ be a
finite-index subgroup of $\Gamma$ so that its transversal is a finite group $K$, making $Q=KP$
the fundamental polyhedron of $G$.
If a path from $Q$ to $gQ$ passes through translates of sides $R_1,\dots,R_m$ of $Q$,
then $g$ can be written as $g=(r_m\dots r_1)(k_1^{-1}\dots k_m^{-1})$, where
$r_i$ is a reflection in a translate of $R_i$ and $k_i\in K$ is the $K$-part of the side-pairing of
$R_i$.
\qed
\end{proposition}

In this paper we work with presentations of groups and their subgroups.
If a group $G$ is given by  presentation $\left<Y\,|\, R\right>$, and $H$ is a finite-index subgroup of
$G$ with Schreier transversal~$X$ (i.e. $G=\cup_{x\in X} Hx$ and every initial segment of an
element $x\in X$ is also an element of $X$),
then the presentation of $H$ can be derived from the presentation of $G$.
Let $g\mapsto \overline g$ be the map $G\to X$ sending $g$ to its coset representative, that is
$Hg=H\overline g$, and let $\gamma:G\to H$ be the map $\gamma(g)=g{\overline g}^{-1}$.

\begin{proposition}[The Reidemeister-Schreier method, \cite{Lyndon-Schupp} Proposition 4.1]
\label{rs}
The presentation of~$H$ is given by $\left<Z\, | \, S\right>$, where
$Z=\{\gamma(xy)\ne 1\ | \ x\in X, y\in Y\}$ and $S$ is the set of all words $xrx^{-1}$ written in terms
of elements of $Z$, for all $x\in X$ and $r\in R$.
\qed
\end{proposition}

\section{The 24-cell, its side-pairings and cycle relations}
\label{24cell}

Working in the ball model of  $\hiv$, let $S_{****}$ be one of the 24 spheres of radius 1 
centered at a point in $\riv$ whose coordinates have two zeroes and two $\pm 1$'s.
The string $*$$*$$*$$*$ consists of symbols $+$,$-$,0 and identifes the sphere by its center.
For example, $S_{+0-0}$ is the sphere centered at $(1,0,-1,0)$.
These spheres are orthogonal to $\bd\hiv$, hence they determine hyperplanes in $\hiv$.
The 24-sided polyhedron $Q$ is defined to be the intersection of the half-spaces corresponding to
those hyperplanes that contain the origin.
The side of $Q$ lying on the sphere $S_{****}$ is also denoted by $S_{****}$.
All the dihedral angles of $Q$ are $\pi/2$ and its 24 ideal vertices are
$v_{\pm000}=(\pm 1,0,0,0)$, $v_{0\pm00}=(0,\pm 1,0,0)$,
$v_{00\pm0}=(0,0,\pm 1,0)$, $v_{000\pm}=(0,0,0,\pm 1)$ and
$v_{\pm\pm\pm\pm}=(\pm 1/2,\pm 1/2,\pm 1/2,\pm 1/2)$.
The polyhedron $Q$ has no real vertices.

Note that two sides intersect if and only if their identifying strings have equal nonzero
entries in one position and the positions of the remaining nonzero entries are different.
Two sides touch at $\bd\hiv$ when they have equal nonzero entries in one position and 
opposite nonzero entries in another position, or if the positions where they have 
nonzero entries are complementary.
Furthermore, an ideal vertex is on a side if its Euclidean distance from the center of the
sphere defining the side is 1.
In the case of vertices whose label has one nonzero entry, a vertex is on a side if and only if the side
has an equal nonzero entry in the same position as the vertex; in the case of vertices
$v_{\pm\pm\pm\pm}$, a vertex  is on a side if and only if the nonzero entries in the label of the side
coincide with the entries in the label of the vertex in the same position.
For example, $S_{0-+0}$ and $S_{+-00}$ intersect, $S_{0-+0}$ and $S_{0+0+}$ are disjoint,
the side $S_{0+0-}$ has six ideal vertices $v_{0+00}$, $v_{000-}$, $v_{\pm+\pm-}$;
$S_{0-+0}$ and $S_{0--0}$ touch at $v_{0-00}$, and $S_{0-+0}$ and $S_{-00-}$ touch
at $v_{--+-}$.

As mentioned in \S\ref{prelims}, every Ratcliffe-Tschantz example has side-pairing isometries
of form~$rk$, where $k\in\z_2^4$ is a composite of  reflections in the coordinate planes
\mbox{$x_1=0,\dots,x_4=0$}.
The label $\pm$$\pm$$\pm$$\pm$ identifies an element of $K$, for example, $k_{+--+}$ is
the composite of reflections in $x_2$ and $x_3$.
Ratcliffe and Tschantz have found in \cite{Ratcliffe-Tschantz} that all sides whose labels have
nonzero entries in the same two positions must have the same $K$-part for their side-pairing.

\rk{Side-pairing labeling convention}
Group letters $a,\dots,l$ and the sides of $Q$ as follows:
\begin{center}
\begin{tabular}{c}
$\{a,b,S_{\pm\pm00}\}$, $\{c,d,S_{\pm0\pm0}\}$, $\{e,f,S_{0\pm\pm0}\}$,\\
$\{g,h,S_{\pm00\pm}\}$, $\{i,j,S_{0\pm0\pm}\}$, $\{k,l,S_{00\pm\pm}\}$.
\end{tabular}
\end{center}
The two letters in each group denote the side-pairings of the four sides from the group.
The first letter always pairs the side labeled $++$, and
the second letter pairs the side whose label is the next unused one, where the symbols
$\pm\pm$ have been ordered in the dictionary order: $++$, $+-$, $-+$, $--$.

The encoding of the side-pairing transformation is done like in \cite{Ratcliffe-Tschantz}.
A string of six characters from $\{1,\dots,9, A,\dots,F\}$, one for each group above, encodes the
$K$-part of each side-pairing.
Each character stands for a hexadecimal number, which is turned into binary form and the
{\it order of digits is reversed}.
Converting every 0 into a plus and every 1 into a minus obtains the label of $k$.
For example, $1$ stands for $k_{-+++}$; $C=12$ stands for $k_{++--}$.
(See Example~\ref{m56}.)
\medskip

The presentation of a group $G$ generated by side-pairings of a polyhedron is obtained by
applying Poincar\'e's polyhedron theorem (\cite{Ratcliffe} Theorem 11.2.2).
Finding presentations of the groups and identifying usable translations in parabolic groups
is crucial to our proof of Theorem~\ref{mainthm}.
An efficient method for these two tasks is essential if one is to deal with a dozen examples,
so we elaborate on it.

\begin{figure}[ht!]\anchor{cube}
\begin{center}
\resizebox{2.5in}{!}{\includegraphics{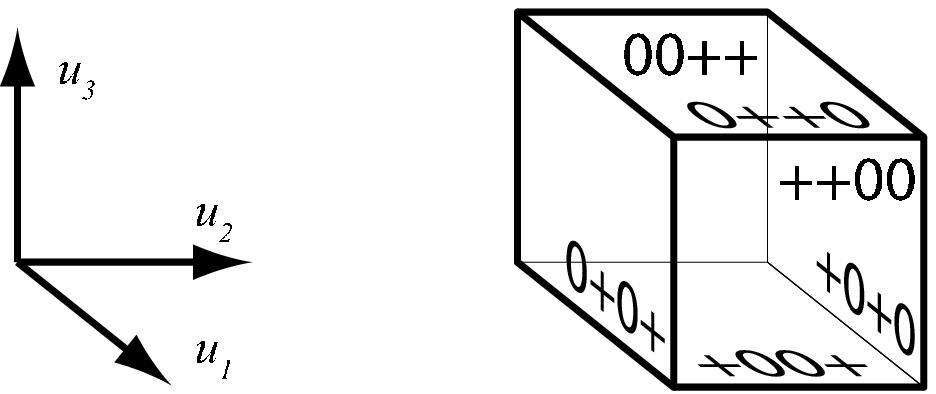}}
\caption{The cube $Q_C$ when $C$ is centered at $v_{++++}$}
\label{cube}
\end{center}
\end{figure}

Let $C$ be a small enough horosphere centered at $v_{++++}$ that intersects only the
sides of~$Q$ that contain $v_{++++}$.
Then $Q_C=C\cap Q$ is a cube with opposing sides \mbox{$S_{++00}$ and $S_{00++}$},
\mbox{$S_{+0+0}$ and $S_{0+0+}$}, $S_{0++0}$ and $S_{+00+}$ (see \figref{cube}).
The tiling of $\hiv$ by translates of~$Q$ intersected with $C$ gives a tiling of $C$ by cubes,
one of which is $Q_C$.

\begin{figure}[ht!]\anchor{ridgecycles}
\begin{center}
\resizebox{4.75in}{!}{\includegraphics{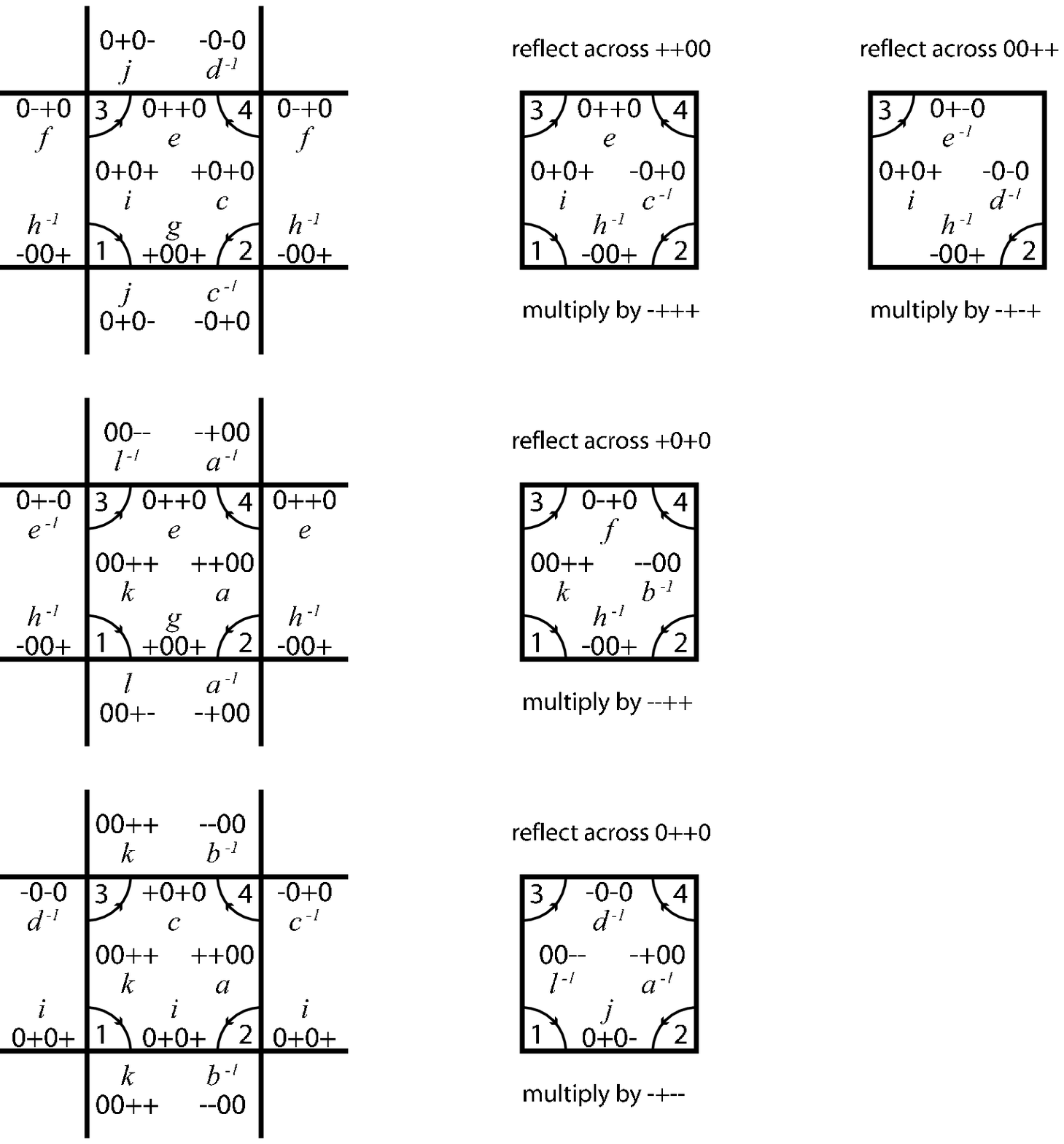}}
\caption{Diagrams that produce cycle relations for $M_{56}$}
\label{ridgecycles}
\end{center}
\end{figure}

We use the term {\it ridge} to denote a codimension-2 face of a polyhedron, while an {\it edge}
is a 1-face.
Every edge of the cube $Q_C$ represents a ridge of $Q$.
Circling around the edge once gives a cycle relation for the presentation in Poincar\'e's theorem.
To easily determine which sides we pass through, use \figref{ridgecycles}.
The left column depicts the intersections of the tiling of $C$ with the three planes that pass
through the center of $Q_C$ and are parallel to its sides.
The side-pairings in that figure are the ones for manifold $M_{56}$ from Example~\ref{m56}.

Let $s$ be the side-pairing sending $S$ to $S'$, $r_S$ the reflection in $S$ and let $k_S$ be
the $K$-part of $s$.
Consider a ridge $U\cap S$, where $U$ is another side of $Q$.
The ridge $U\cap  S$ is sent via $s$ to a ridge $V\cap S' $, where $V$ is another side of~$Q$.
Now $V\cap S'=k_S r_S (U\cap S) =k_S (U\cap S)=k_S(U)\cap S'$ so we conclude $V=k_S(U)$.
Therefore, adjacent to $Q$ on the other side of $S$ is $s^{-1}Q$, whose sides $s^{-1}V$
and $s^{-1}S'$ meet at ridge $U\cap S=s^{-1}(V\cap S')$.

The observation $V=k_S(U)$ allows us to easily find the labelings on all the struts in the diagrams
in the left column of \figref{ridgecycles} (labelings inside the square come from $Q_C$).
For example,  consider the sides around the vertex labeled 1 in the topmost diagram.
The bottom left strut is the translate of the side
$k_{S_{+00+}}(S_{0+0+})=k_{-++-}(S_{0+0+})=S_{0+0-}$, whose side-pairing is $j$.
The left bottom strut is the translate of $k_{S_{0+0+}}(S_{+00+})=k_{--++}(S_{+00+})=S_{-00+}$,
whose side-pairing is $h^{-1}$.

Cycle relations are words obtained circling an edge of $Q_C$ and successively adding the
side-pairings at the end of the word.
Care must be taken to add $s^{-1}$ at the end of the word whenever we {\it exit} through 
a side labeled $s$ and to add $s$ whenever we {\it enter} into a side labeled~$s$.
Going in the direction of the arrow around vertex 1 we exit through $g$ and $j$ and then enter
into $h^{-1}$ and $i$, so the cycle relation is $g^{-1}j^{-1}h^{-1}i=1$.

We thus obtain 12 relations for $G$ from the first column of \figref{ridgecycles}.
Because $Q$ has 96 ridges that fall into 24 ridge cycles, each corresponding to a relation, 12 further
relations are needed.
They are found by going around other edges in the tiling of $C$.
If we reflect the three planes we used for sections of $Q_C$ in the sides of $Q_C$ parallel to them,
we get new sections of the tiling represented by diagrams in the middle and right columns
of \figref{ridgecycles}.
The fact $V=k_S(U)$ from above is now used to see that the labelings on the new diagrams are
obtained by applying $k_S$ to every label in the original diagram, where $S$ is the side used to
reflect the section plane.

The middle and right columns do not include struts because cycle relations follow
from only the labels on the square part and the cycle relations for the left column.
If $A=\{S_1,kS_1,S_2,kS_2\}$ is the set of four sides whose labels all have the same nonzero
positions, and are hence paired by the same element $k\in K$, we note that any element $k'\in K$ 
acts on $A$ and, due to commutativity of $K$, preserves or interchanges the subsets of $A$
containing paired sides.
Along with the fact that $k'^2=1$ for every $k'\in K$ this means that action of $k'$ on $A$ can be
inferred just from knowing $k'S_1$.
For example, in the third row in Fig. 1, we see that by applying $k_{S_{0++0}}$ to the sides labeled
$i$ and $a$ we get sides labeled $j$ and~$a^{-1}$.
We infer that applying $k_{S_{0++0}}$ to sides labeled $j$ and $b$ gives sides labeled $i$ and $b^{-1}$.
Thus, the cycle relation for vertex 2 in the new diagram is the cycle relation $i^{-1}b^{-1}ia$ for vertex~2
in the old diagram with letters replaced in the above pattern, obtaining $j^{-1}bja^{-1}$.
Now, if  a cycle relation coming from the new diagram is the inverse, a cyclic permutation
or a combination of the two of a corresponding cycle relation in the old diagram, we have found
the same ridge cycle, and have to reflect the section plane in the other side of the cube $Q_C$.

We argue that the nine diagrams thus obtained will contain all the cycle relations (multiple times, in
fact).
Note that in $Q$ every ridge is an ideal triangle, two of whose ideal vertices are of type
$v_{\pm\pm\pm\pm}$.
For example,  the ridge $R=S_{0+0+}\cap S_{0++0}$ is an ideal triangle with vertices
$v_{0+00}$, $v_{++++}$ and $v_{-+++}$ and edges
$E_1=S_{0+0+}\cap S_{0++0}\cap S_{00++}$,
$E_2=S_{0+0+}\cap S_{0++0}\cap S_{++00}$, and
$E_3=S_{0+0+}\cap S_{0++0}\cap S_{-+00}$.
The intersection of translates of edges of $Q$ with $C$ are the vertices in the tiling of $C$.
An edge of the cube $Q_C$ is an arc in the ideal triangle going from one  edge of the triangle to
another.
Because elements of $K$ do not change the nonzero positions in the labeling of sides, 
only $E_2$ and $E_3$ might be in the same cycle of edges.
Thus, in the manifold $\hiv/G$ the ideal triangle might only have its edges $E_2$ and $E_3$ 
identified, since no points on $E_1$ are identified.
If  ridge $R$ and the ridge represented by reflection of $R\cap C$ in the side $S_{00++}$ of $Q_C$
were in the same cycle, the segment consisting of $R\cap C$ and its reflection would
project to a smooth curve in the ideal triangle that crosses $E_1$ twice, which is impossible.
Since in $Q$ there are only 8 ridges of type $S_{0\pm0\pm}\cap S_{0\pm\pm0}$ we conclude that
every one of them is either in the cycle of $R$, or the cycle of the ridge represented by the reflection
of $C\cap R$ in $S_{00++}$.
Reasoning similarly for all the other ridge cycles we see that the opportunity for cycle relation
repetitions occurs only for ridges represented by vertices labeled 2 and 3 in the middle column.
Often there are no repetitions: in Example~\ref{m56} we needed three section planes only in the
first row.

Now we describe how to identify and choose suitable translations in parabolic subgroups.
Let $G_v$ be the parabolic group that is the stabilizer of an ideal vertex $v$ of $Q$, and let
$C$ be a horosphere centered at $v$.
The intersection $Q_C=Q\cap C$ is always a Euclidean cube.
Position a coordinate system with origin in the center of $Q_C$ and axes
perpendicular to its sides, and consider the tiling of $C$ by translates of $Q_C$.
To simplify language, we identify a cube~$gQ_C$ and the element $g$.
According to Proposition~\ref{moving}, by considering a path from $Q_C$ to~$gQ_C$ that passes
only through sides of cubes (e.g. a path that is piecewise parallel to the coordinate axes), we may
write $g=(r_m\dots r_1)(k_1^{-1}\dots k_m^{-1})$.
Since $r_m\dots r_1$ preserves~$v$, $g\in G_v$ if and only if $k_1^{-1}\dots k_m^{-1}v=v$.

\begin{figure}[p]\anchor{trans1-4}
\begin{center}
\resizebox{4.8in}{!}{\includegraphics{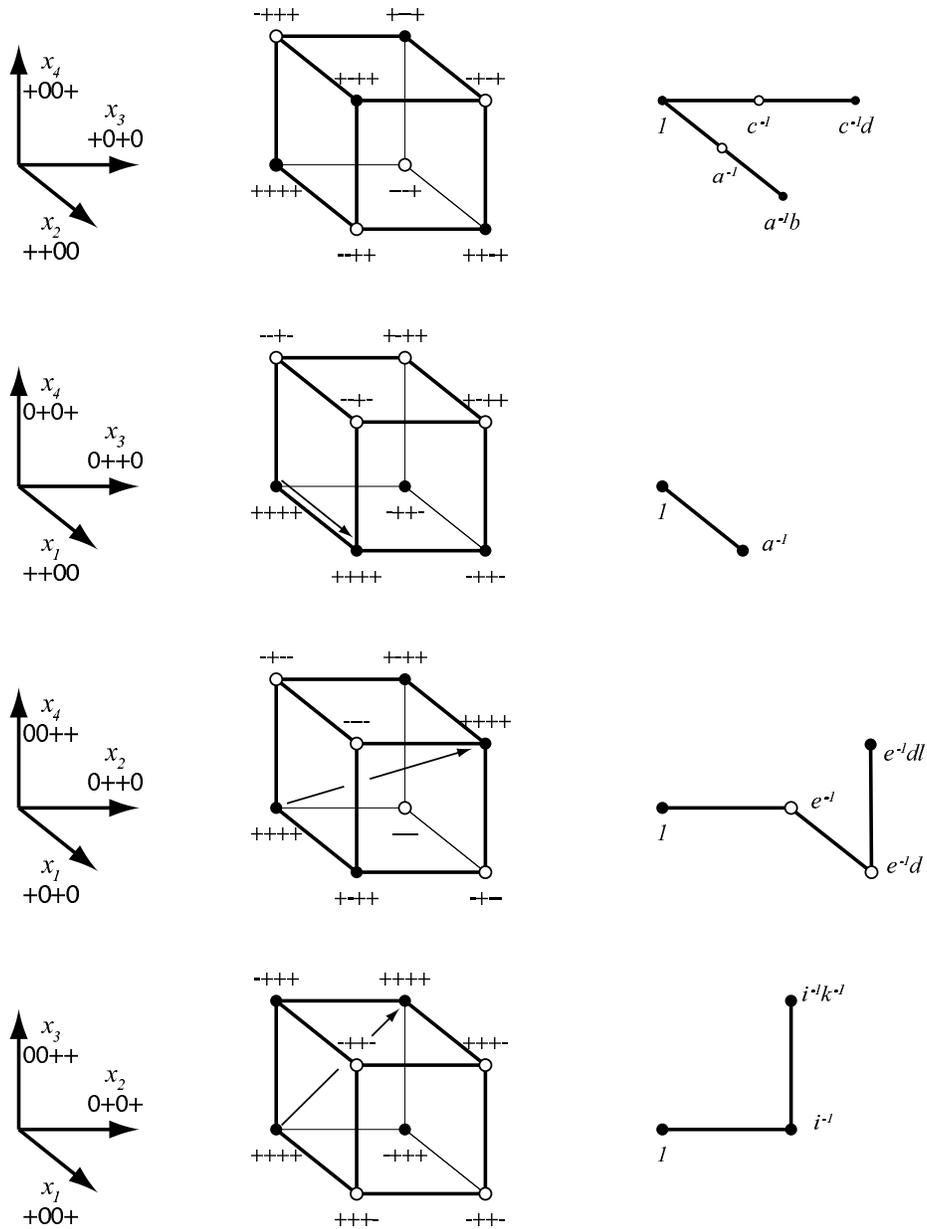}}
\caption{Translations in $8Q_C$ for manifold $M_{56}$}
\label{trans1-4}
\end{center}
\end{figure}

Let $v=v_{+000}$, $v_{0+00}$, $v_{00+0}$ or $v_{000+}$.
The coordinate axes perpendicular to sides of $Q_C$ may be identified with coordinate axes
$x_1$, $x_2$, $x_3$ and $x_4$, in the arrangement suggested by \figref{trans1-4}.
The 4-symbol labels on the axes indicate which side we pass through if we move from the center
of $Q_C$ in the direction of the arrow.
Parallel sides of each $Q_C$ have same nonzero positions but differ in sign in positions other
than the nonzero position of the vertex label (a side and a vertex have the same symbol
in the nonzero position of the vertex label when the vertex is on the side).
Moving two steps in the directions of the axes produces translations since $k_1=k_2$ in
every case.
Let $T_v$ be the subgroup of $G_v$ generated by those three translations and let $TG_v$ be the
translation subgroup of $G_v$.
For our purposes we look for translations only in $8Q_C$, the $2\times 2 \times 2$ block of cubes 
that is the fundamental polyhedron for~$T_v$.
A Euclidean transformation $Ax+a$ is a translation if and only if its rotational part~$A$ is the identity.
We get the rotational part of $r_m\dots r_1 k_1^{-1}\dots k_m^{-1}$ by replacing $r_i$ by an
element~$k_{r_i}$ that represents the reflection in the coordinate plane parallel to the reflection plane
of~$r_i$ and computing the product.
Note that it is enough to find the rotational parts of  the three cubes adjacent to $Q_C$ in the
direction of the three axes.
Any other rotational part is obtained by multiplying those three: every time we move
by one step in the direction of an axis, we multiply by the rotational part of one of the three
cubes that corresponds to that axis.
The dots in \figref{trans1-4} represent centers of cubes comprising~$8Q_C$;
the dot is filled if it corresponds to an element preserving~$v$.
Labels on the dots give the rotational part of the group element corresponding to the dot.
The three zero positions in the vertex label give the rotational part of the transformation, while the
nonzero position tells whether the transformation preserves $v$.
For example, in the topmost diagram in \figref{trans1-4} the label $+$$+$$-$$+$ on a vertex reveals
an element that preserves $v$ whose rotational part is a reflection in the $x_3$-plane.
The label $-$$-$$+$$+$ reveals that the corresponding element does not preserve~$v$.
Arrows in the middle column indicate translations found in $8Q_C$, they correspond to dots
labeled $+$$+$$+$$+$.
The right column of the diagram computes the group elements corresponding to the dots
of interest.
This is done by considering which translates of sides we pass through as we go from~$Q_C$
to the dot.

\begin{figure}[ht!]\anchor{trans5}
\begin{center}
\resizebox{4.5in}{!}{\includegraphics{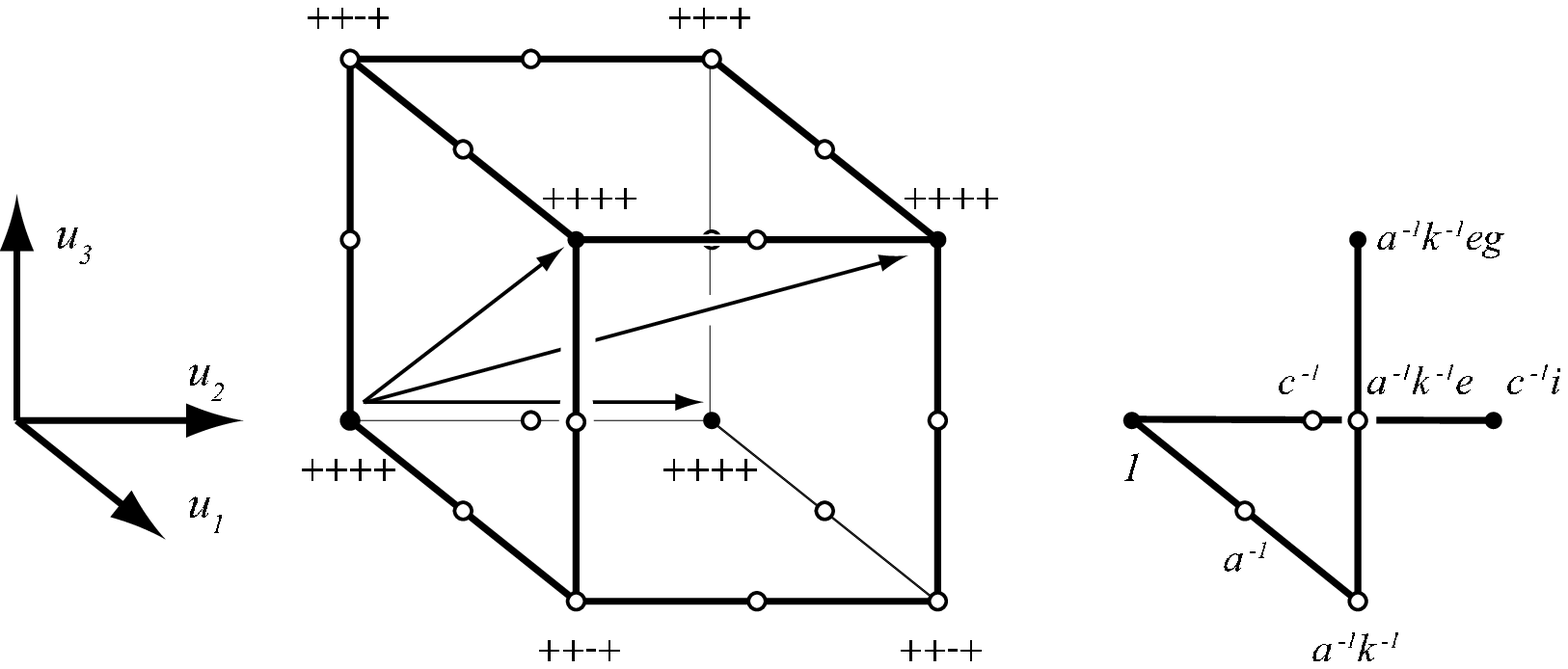}}
\caption{Translations in $64Q_C$ for manifold $M_{56}$}
\label{trans5}
\end{center}
\end{figure}

Let $v=v_{++++}$, then $Q_C$ is the cube from \figref{cube}.
Here the condition $k_1^{-1}\dots k_m^{-1}v=v$ becomes $k_1^{-1}\dots k_m^{-1}=1$, since
$K$ acts freely on the set of vertices $v_{\pm\pm\pm\pm}$.
Moving four steps in the directions of the three axes produces translations because
$k_1=k_3$, $k_2=k_4$ and $r_1\dots r_4$ is
a composite of reflections in four parallel planes.
Let $T_v$ be the subgroup of~$G_v$ generated by those three translations and $TG_v$ be the
translation subgroup of $G_v$.
For our purposes, we look for translations in $64Q_C$, the $4\times 4 \times 4$ block of cubes
that is the fundamental polyhedron for $T_v$.
It is clear that we get a translation if and only if we move an even number of steps in the direction
of each axis, and have $k_1^{-1}\dots k_m^{-1}=1$.
Therefore, we have to check $k_1^{-1}\dots k_m^{-1}=1$ for the 7 cubes obtained by moving $Q_C$
by a combination of two units in every direction.
The dots in \figref{trans5} depict the centers of some cubes in~$64Q_C$, and the ones
at the vertices of the cube are the interesting ones.
Labels indicate the products of the $K$-parts, so  a $+$$+$$+$$+$ label means the dot
corresponds to a translation.
At right we compute the group elements corresponding to the dots of interest.
One of the translations is omitted, since it is the sum of the other two in the diagram.

\section{Examples of ``link'' complements in the 4-sphere}
\label{examples}

In this section we show that 12 nonorientable Ratcliffe-Tschantz manifolds have orientable double
covers that are complements of tori and Klein bottles in $S^4$.
First of all, we need criteria that rule out a manifold from having this property.
The first criterion is:

\begin{proposition}
\label{crithom}
Let $M$ be a complement of $r$ tori and $s$ Klein bottles in $S^4$.
Then (coefficients are in $\z$):
$$
H_1 M = \z^r \oplus \z_2^s,\ 
H_2 M =\z^{2r} \oplus \z^s,\ 
H_3 M= \z^{r+s-1}.
$$
\end{proposition}

\begin{proof}
The proof is a simple application of Alexander's duality: if $A$ is a compact submanifold
of $S^n$, then $\tilde H_i(S^4-A)\cong \tilde H^{3-i}(A)$.
\end{proof}

The other criterion was developed in \cite{Ivansic3} and exploits a homomorphism that exists for
every group $G$ of a Ratcliffe-Tschantz manifold.
It is defined by $\phi:G\to \z_2^6$, $a,b\mapsto e_1$, $c,d\mapsto e_2,\dots, k,l\mapsto e_6$,
where $e_i$ is a canonical generator of $\z_2^6$.

\begin{proposition}
\label{critphi}
Let $T$ be the set of all translations in $G$.
If $\dim_{\z_2}\left<\phi(T)\right> <5$ then the orientable double cover of $M$ is not a complement
in the 4-sphere.
\qed
\end{proposition}

Using a computer to test the 1149 nonorientable Ratcliffe-Tschantz manifolds for both criteria,
we found that only 49 remained eligible to have double-cover complements in $S^4$, namely,
manifolds numbered
\begin{displaymath}
\begin{array}{rrrrrrrrrrrrr}
23& 25 & 28& 29& 32& 33& 35& 36& 40& 46& 47& 51& 56 \\
60& 71 & 72 & 73 & 74& 80& 92& 93 & 94& 96& 109 & 112 &116\\
121& 131& 132& 139& 170& 171& 231& 235& 236& 251& 292& 296& 297\\
299&426& 427& 434& 435& 1011& 1091& 1092& 1094& 1095 &&&
\end{array}
\end{displaymath}
If the double cover of any of these manifolds is a complement in $S^4$, the potential link
structure is visible from its homology by Proposition~\ref{crithom}.
The computer calculation showed that only twelve different link structures are possible.
We then proceeded to prove (manually) that for each of those link structures there is at least
one manifold from the above list whose orientable double cover is the complement of such a
link in $S^4$.
Those manifolds are listed in Theorem~\ref{mainthm}.
Most likely, other manifolds from the list above have the same property, but it is not  easy
to ascertain that we will be getting examples different from the ones in Theorem~\ref{mainthm},
since orientable double covers of nonisometric manifolds may be isometric.

Let $G$ be any group, $H$ a normal subgroup of $G$ and let $\ncl A\ncr_H$ be the normal closure
of a set $A\subset H$ in $H$.
If $G$ is the biggest group under consideration, $\ncl A \ncr$ or $\ncl A \ncr_G$ will denote the
normal closure of $A$ in $G$.
Clearly, $\ncl A\ncr_H$ is the set of all elements of the form $\prod h_i a_i h_i^{-1}$, where $h_i\in H$,
$a_i\in A\cup A^{-1}$.
If $X$ is a transversal of $H$ in $G$ it is easy to see that for a set $A\subset H$,
$\ncl A \ncr_G=\ncl xAx^{-1},\ x\in X\ncr_H$.
Note that $\ncl A \ncr_G$ is a subgroup of $H$, since $H$ is normal in $G$.

While $\ncl A \ncr_H\subseteq \ncl A \ncr_G$, in general, $\ncl A \ncr_H \ne \ncl A \ncr_G$.
For example, if $H=\left<a\right>*\left< b\right>$, $G=H\rtimes\z_2$, where $\z_2$ acts on $H$
by swapping the generators $a$ and $b$, then clearly $\ncl a \ncr_G=H$, while $\ncl a \ncr_H$
is a subgroup of $H$ with infinite index.

The following proposition will prove useful.

\begin{proposition}
\label{normclosures}
Let $H$ be a normal subgroup of $G$ with transversal $X$.
Suppose a set $A$ has the property that $\ncl A\ncr_H=\ncl A\ncr_G$.
Let $B\subset H$ have the property: for every $b\in B$ and every $x\in X$, $[xbx^{-1}]=[b^{\pm1}]$,
where $[\ ]$ denotes classes in $G/\ncl A \ncr_G$.

Then $\ncl A\cup B \ncr_H=\ncl A\cup B \ncr_G$.
\end{proposition}

\begin{proof}
Elements of type $hxa(hx)^{-1}$ and $hxb(hx)^{-1}$ generate the subgroup $\ncl A\cup B\ncr_G$,
where $a\in A\cup A^{-1}$, $b\in B\cup B^{-1}$, $x\in X$, $h\in H$.
Now $\ncl xAx^{-1},\ x\in X \ncr_H=\ncl A \ncr_G=\ncl A \ncr_H$ implies
that every element  $xax^{-1}$ is a product of elements of type $hah^{-1}$, so every
$hxa(hx)^{-1}=h(xax^{-1})h^{-1}$ has the same property.
Since $[xbx^{-1}]=[b^{\pm1}]$ in $G/\ncl A \ncr_G$ means that $xbx^{-1}=b^{\pm1}u$, for some
$u\in \ncl A \ncr_G=\ncl A \ncr_H$,  $xbx^{-1}$ is a product of $b^{\pm1}$ and elements
of type $hah^{-1}$, making $hxb(hx)^{-1}$ a product of elements of type $hbh^{-1}$ and $hah^{-1}$.
Therefore, $\ncl A\cup B\ncr_G$ is generated by elements of type $hah^{-1}$ and $hbh^{-1}$,
so $\ncl A\cup B\ncr_G=\ncl A\cup B\ncr_H$.
\end{proof}

\begin{theorem}
\label{mainthm}
The orientable double covers of the nonorientable Ratcliffe-Tschantz
manifolds listed in \tableref{table1} are complements of the indicated
combination of tori and Klein bottles in a manifold that is
homeomorphic to the 4-sphere.\end{theorem}

\begin{table}\small\anchor{table1}
\begin{center}
\begin{tabular}{l|l}
The orientable double & is a complement \\
cover of manifold &   in $S^4$ of:\\
\hline
no.~1011 & \hskip5pt 5 tori \\
no.~71 & \hskip5pt 6 tori\\
no.~23 & \hskip5pt 7 tori\\
no.~1092  & \hskip5pt 8 tori\\
no.~1091 & \hskip5pt 9 tori \\
no.~231 & \hskip5pt 3 tori, 4 Klein bottles \\
no.~112&  \hskip5pt 4 tori, 2 Klein bottles \\
no.~56& \hskip5pt 4 tori, 3 Klein bottles \\
no.~92 & \hskip5pt 5 tori, 1 Klein bottle\\
no.~51 & \hskip5pt 5 tori, 2 Klein bottles \\
no.~40& \hskip5pt 6 tori, 1 Klein bottle \\
no.~36 & \hskip5pt 6 tori, 2 Klein bottles\\
\end{tabular}
\end{center}
\nocolon\caption{}\label{table1}
\end{table}

\begin{remark}
Note that  Theorem~\ref{mainthm} claims that double covers of various Ratcliffe-Tschantz
manifolds are complements in a manifold that is {\it homeomorphic} to the 4-sphere, rather
than diffeomorphic.
This is because we used Freedman's theory, applicable only to the topological category,
to recognize a 4-sphere.
It is still unknown whether manifolds homeomorphic to the 4-sphere are diffeomorphic to the $S^4$
with the standard differentiable structure.
Using Kirby diagrams of a 4-manifold, Ivan\v si\'c has shown that the topological $S^4$ inside
which $\tilde M_{1011}$ is a complement of 5 tori is indeed diffeomorphic to the standard
differentiable $S^4$.
We anticipate that the same is true for the examples treated here.
\end{remark}

\begin{proof}
Let $M$ be a nonorientable Ratcliffe-Tschantz manifold, $\tilde M$ its orientable double cover,
$G=\pi_1 M$, $H=\pi_1 \tilde M$.
Clearly, $H$ is the index-2 subgroup of $G$ containing all the orientation preserving isometries.
If we can find primitive normal translations $t_1,\dots,t_m$, one in each boundary component of $\tilde M$,
then $\tilde M$ is a complement in a closed manifold $N$.
As was shown in \cite{Ivansic2}, the Euler characteristics of $\tilde M$ and $N$ are equal.
By Proposition~\ref{presentationofN}  $\pi_1 N=H/\ncl t_1,\dots,t_m\ncr_H$.
The classification of 4-dimensional simply-connected manifolds (see \cite{Freedman-Quinn} or
\cite{Gompf-Stipsicz}) gives that if $\pi_1 N=1$ and $\chi(N)=2$, then $N$ is homeomorphic to $S^4$.
Since $\chi(\tilde M)=2$, it will be enough to show that $H/\ncl t_1,\dots,t_m\ncr_H=1$.

The presentation of $G$ has 12 generators and 24 relations so the presentation of the \mbox{index-2}
subgroup $H$ obtained by the Reidemeister-Schreier method will have 23 generators and
48 relations,  very tedious to compute with.
Thus, the challenge is to show that $H/\ncl t_1,\dots,t_m\ncr_H=1$ while working with the
presentation of $G$ as much as possible.

The ten homeomorphism types of flat 3-manifolds are identified by the letters A--J, where the
manifolds are ordered in the same way as in \cite{Hantzsche-Wendt} or \cite{Wolf}.
Every nonorientable Ratcliffe-Tschantz manifold that passes both criteria from above has boundary
components of type A, B, G, H, I or J.
Type A is the 3-torus and type B is an (orientable) $S^1$-bundle over a Klein bottle.
The nonorientable types G and H have type A as their orientable double cover and the nonorientable
types I and J have type B as their orientable double cover.
Holonomy groups of boundary components of $M$ are computed in the course of searching
for suitable translations.
They help determine whether chosen translations are normal in a boundary component
and they distinguish types A, B, G or H, and I or J.
Although not essential for our proof, types G and H and can further be distinguished by the fact that
at least one generator of the translation subgroup of H is not a normal translation.
There is a criterion that distinguishes I and J as well.

Note that orientable boundary components of $M$ lift to two homeomorphic boundary components
in $\tilde M$ and nonorientable ones lift to their orientable double covers.
Thus, knowing the boundary components of $M$ immediately gives the link structure that
$\tilde M$ is a complement of (a torus or a Klein bottle for every component of type A or B,
respectively).

Translations in boundary components of $\tilde M$ have to be carefully chosen in order
to be $S^1$-fibers, to produce $H/\ncl t_1,\dots,t_m\ncr_H=1$ and to allow us to work in $G$.
It was shown in \cite{Ivansic3} that the homomorphism $\phi$ from Proposition~\ref{critphi}
sends to 0 every subgroup $T_v$ in~\S\ref{24cell}, and that $H/\ncl t_1,\dots,t_m\ncr_H=1$
implies $\dim_{\z_2}\left<\phi(t_1),\dots,\phi(t_m)\right>\ge 5$.
Therefore, one should choose translations in $8Q_C$ or $64Q_C$ whenever possible, since one
risks otherwise that $\dim_{Z_2} \left<\phi(t_1),\dots,\phi(t_m)\right>$ comes short of 5.
A translation of type $t=t'+t''$  where $t'\in T_v$ and $t''\in 8Q_C$ or $64Q_C$, is a possible choice,
but we avoid them, because its length in terms of generators $a$--$l$ is big, possibly
complicating the computation.

Let $E$ be a boundary component of $M$ and $\tilde E$ a component of its lift in $\tilde M$.
In order for $\tilde M$ to be a complement, we must choose a translation $t$ that is normal
and primitive in $\tilde E$.
Note that $t$ need not be normal in $E$, but it is very useful if it is when $E$ is nonorientable.
For in that case, there is an orientation-reversing $x\in G$ so that $xtx^{-1}=t^{\pm 1}$, which by
Proposition~\ref{normclosures} implies $\ncl A, t\ncr_H=\ncl A, t \ncr_G$ for any subset
$A$ of $H$ satisfying $\ncl A \ncr_H=\ncl A\ncr_G$.
This allows us to continue computing in $G$.
We introduce a term for this property: we say that $t$ is {\it self-conjugate} if there
exists an orientation-reversing $x\in G$ so that $xtx^{-1}=t^{\pm 1}$.

Checking that $t$ is primitive in $\tilde E$ is easy because $\tilde E$ is either of type A or B.
In the first case, $t$ must be the shortest translation in its direction, in the second, it must not
point in the same direction as the axis of rotation by $\pi$ that generates the holonomy group.

If $E$ is orientable, its lift has two components $\tilde E$ and $\tilde E'$, whose fundamental groups
are conjugate by an orientation reversing $x\in G$.
Choosing a primitive normal translation $t$ in $\tilde E$ gives the obvious choice $xtx^{-1}$ for
$\tilde E'$.
This choice works most of the time, but on occasion, like in Example~\ref{m36}, one must
choose a different translation in $\tilde E'$ in order to satisfy the condition
$\dim_{\z_2}\left<\phi(t_1),\dots,\phi(t_m)\right>\ge 5$.

To prove the theorem, each example needs to be handled separately.
We show three of the most difficult examples in decreasing detail and indicate how to proceed
on others.

\begin{example}
\label{m56}
The side-pairing of manifold $M_{56}$ is encoded by 13D935.
According to our side-pairing convention, the side-pairings are (the $K$-parts are under the arrow):
\begin{displaymath}
\begin{array}{ll}
S_{++00} \arrtop{a}{-+++} S_{-+00}\hskip10pt
&
S_{+-00} \arrtop{b}{-+++} S_{--00}\hskip10pt
\\
\rule[18pt]{0pt}{0pt}%
S_{+0+0} \arrtop{c}{--++} S_{-0+0}\hskip10pt
&
S_{+0-0} \arrtop{d}{--++} S_{-0-0}
\\
\rule[18pt]{0pt}{0pt}%
S_{0++0} \arrtop{e}{-+--} S_{0+-0}\hskip10pt
&
S_{0-+0} \arrtop{f}{-+--} S_{0--0}\hskip10pt
\\
\rule[18pt]{0pt}{0pt}%
S_{+00+} \arrtop{g}{-++-} S_{-00-}\hskip10pt
&
S_{+00-} \arrtop{h}{-++-} S_{-00+}
\\
\rule[18pt]{0pt}{0pt}%
S_{0+0+} \arrtop{i}{--++} S_{0-0+}\hskip10pt
&
S_{0+0-} \arrtop{j}{--++} S_{0-0-} \hskip10pt
\\
\rule[18pt]{0pt}{0pt}%
S_{00++} \arrtop{k}{-+-+} S_{00-+}\hskip10pt
&
S_{00+-} \arrtop{l}{-+-+} S_{00--}
\end{array}
\end{displaymath}
The orientation-reversing generators are $c$, $d$ and $g$--$l$ (those whose $K$-part has an even
number of minuses in its label).

We use diagrams in \figref{ridgecycles} to find the cycle relations as described in \S\ref{24cell}.
The arrangement of groups of relators below corresponds to the arrangement of diagrams
used to obtain them and the numbers 1--4 are labels of the vertices in the diagrams corresponding
to each cycle relation.
The cycle relations marked by $\bullet$ are repetitions of those from the left column, so
we reflect the section plane in $S_{00++}$ and obtain the two new relations in the right
column.
\begin{displaymath}
\begin{array}{rrr}

\begin{array}{rr}
1)& g^{-1}j^{-1}h^{-1}i=1 \\
2)& g^{-1}ch^{-1}c=1\\
3)& e^{-1}j^{-1}fi= 1\\
4)&  e^{-1}dfc=1          
\end{array}
&
\begin{array}{r}
 hj^{-1}gi=1 \\
\bullet\   hc^{-1}gc^{-1}=1\\
 \bullet\  e^{-1}j^{-1}fi= 1\\
 e^{-1}d^{-1}fc^{-1}=1
\end{array}
&
\begin{array}{r}
\\
hd^{-1}gd^{-1}=1\\
ej^{-1}f^{-1}i=1\\
\ 
\end{array}

\\
&&\\

\begin{array}{rr}
1)& g^{-1}l^{-1}h^{-1}k=1\\
2)& g^{-1}ah^{-1}a=1\\
3)& e^{-1}le^{-1}k=1\\
4)& e^{-1}aea=1
\end{array}
&
\begin{array}{r}
hl^{-1}gk = 1\\
hb^{-1}gb^{-1}=1\\
f^{-1}lf^{-1}k=1\\
f^{-1}b^{-1}fb^{-1}=1
\end{array}
& 

 \\
&&\\

\begin{array}{rr}
1)& i^{-1}k^{-1}ik=1\\
2)& i^{-1}bia=1\\
3)& c^{-1}k^{-1}d^{-1}k =1\\
4)& c^{-1}bc^{-1}a=1
\end{array}
&
\begin{array}{r}
j^{-1}ljl^{-1}=1\\
j^{-1}b^{-1}ja^{-1}=1\\
dlcl^{-1}=1\\
db^{-1}da^{-1}=1
\end{array}
&
\end{array}
\end{displaymath}

The diagrams in Figures~\ref{trans1-4} and~\ref{trans5} find the translations of interest for $M_{56}$.
The following table identifies boundary component types and the translations.
$\tilde M_{56}$ has boundary components of type BBBAAAA, making it a complement
of 3 Klein bottles and 4 tori.
Boundary components $E_1$--$E_4$ are the ones corresponding to ideal vertices
$v_{+000}$--$v_{000+}$, while $E_5$ corresponds to vertex $v_{++++}$.
Elements of the holonomy group are identified by their planes of reflection, for example,
$x_2 x_3$ is the composite of reflections in planes $x_2$ and $x_3$.
\begin{center}
\begin{tabular}{c|c|c|c|c}
comp. & type & translations in $8Q_C$ or $64 Q_C$ & holonomy & self-conj?\\
\hline\rule[15pt]{0pt}{0pt}%
$E_1$ & $I$ & none &  $x_2$, $x_3$, $x_2x_3$& n/a\\
$E_2$ & $B$ & $a$ & $x_1x_4$ & no\\
$E_3$ & $H$ & $e^{-1}dl$ & $x_2$ & no\\
$E_4$ & $G$ & $i^{-1}k^{-1}$ & $x_1$ & yes\\
$E_5$ & $A$ & $c^{-1}i$, $a^{-1}k^{-1}eg$ & trivial & no
\end{tabular}
\end{center}
\figref{trans1-4} shows why $e^{-1}dl$ is not normal in $E_3$: reflecting the vector in the
third cube in the plane $x_2=0$ produces a nonparallel vector.
Similarly, $i^{-1}k^{-1}$ is normal in $E_4$ because it is invariant under reflection in the
plane $x_1=0$.
Note that there are no translations in $E_1$ that are in the cube $8Q_C$, so we
use one of the translations $a^{-1}b$, $c^{-1}d^{-1}$, or $g^{-1}h$ that generate $T_{v_{+000}}$.
The nontrivial element of holonomy of $\tilde E_1$ is a rotation in the $x_2 x_3$-plane.
All three translations are thus normal in $\tilde E_1$, but the third one is not primitive, since
it is a power of a slide-rotation.
Make the following choices for translations:
\begin{displaymath}
\begin{array}{lllll}
t_1=c^{-1}d,  & t_2=a,  &  t_3=e^{-1}dl,   & t_4=i^{-1}k^{-1}, &  t_5=c^{-1}i, \\
         &              t'_2=ct_2c^{-1}, & & & t'_5=c(a^{-1}k^{-1}eg)c^{-1}.
\end{array}
\end{displaymath}
Because $t_1$ and $t_4$ are self-conjugate, we have
$\ncl t_1,t_2,t_4 \ncr_G=\ncl t_1,t_2,t'_2,t_4 \ncr_H$.
Adding $t_1=t_2=t_4=1$ to the relators we get $d=c$, $k=i^{-1}$ and $a=1$, which immediately
yields $b=1$ and  $h=g^{-1}$, so the presentation for $G/\ncl t_1,t_2,t_4 \ncr$ has generators
$c$, $f$, $g$, $i$, $j$, $l$ and relations (repeating or trivial relations omitted):
\begin{displaymath}
\left\{
\begin{array}{r}
g^{-1}j^{-1}gi=1\\
g^{-1}cgc=1\\
e^{-1}j^{-1}fi=1\\
e^{-1}cfc=1\\
ej^{-1}f^{-1}i=1\\
e^{-1}c^{-1}fc^{-1}=1
\end{array}
\right.
\hskip10pt
\left\{
\begin{array}{r}
g^{-1}l^{-1}gi^{-1}=1\\
e^{-1}le^{-1}i^{-1}=1\\
f^{-1}lf^{-1}i^{-1}=1
\end{array}
\right.
\hskip10pt
\left\{
\begin{array}{r}
c^{-1}ic^{-1}i^{-1}=1\\
c^{-2}=1\\
j^{-1}ljl^{-1}=1\\
clcl^{-1}=1
\end{array}
\right.
\end{displaymath}
Taking quotients is equivalent to adding relators so, as is customary in working with group
presentations, we usually omit class notation.
Taking into account $c^2=1$, note that one relation in the group at right says that $i$ and $c$
commute in $G/\ncl t_1,t_2,t_4 \ncr$, and one relation in the left group says $g$ and $c$
commute ($c=c^{-1}$).
But then $ct_5c^{-1}=cc^{-1}ic^{-1}=t_5$, fulfilling conditions of Proposition~\ref{normclosures},
so we conclude $\ncl t_1,t_2,t_4,t_5 \ncr_G=\ncl t_1,t_2,t'_2,t_4,t_5\ncr_H$.

In the group $G/\ncl t_1,t_2,t_4,t_5 \ncr$, we now have $i=c$, so $j=gig^{-1}=gcg^{-1}=c$
and  $l=gi^{-1}g^{-1}=gc^{-1}g^{-1}=c$, which in conjuction with  $e^{-1}le^{-1}i^{-1}=1$
implies $ e^{-1}ce^{-1}c^{-1}=1$, in other words, $cec^{-1}=e^{-1}$.
Since now $t_3=e^{-1}$, we have $ct_3c^{-1}=ce^{-1}c^{-1}=e=t_3^{-1}$, so
Proposition~\ref{normclosures} yields
$\ncl t_1,t_2,t_3,t_4,t_5 \ncr_G$ $=\ncl t_1,t_2,t'_2,t_3,t_4,t_5 \ncr_H$.

The presentation of the group $G/ \ncl t_1,t_2,t_3,t_4,t_5\ncr$ is obtained from the above relations
combined with $i=c$ and $e=1$ (following from $t_5=t_3=1$).
Then \mbox{$f=1$} follows.
Using $j=l=c$, the presentation simplifies to
\begin{displaymath}
\left< c,g\, |\, c^2=1,\ cgc^{-1}=g\right>=\z\oplus\z_2.
\end{displaymath}
In this group, $t'_5=ca^{-1}k^{-1}egc^{-1}=c1c1gc^{-1}=gc^{-1}=gc$, so $ct'_5c^{-1}=t'_5$.
Proposition~\ref{normclosures} says
$\ncl t_1,t_2,t_3,t_4,t_5,t'_5\ncr_G=\ncl t_1,t_2,t'_2,t_3,t_4,t_5,t'_5 \ncr_H$.
Now $G/ \ncl t_1,t_2,t_3,t_4,t_5,t'_5\ncr=\left<c\,|\,c^2=1\right>=\z_2$, so
$H/\ncl t_1,t_2,t'_2,t_3,t_4,t_5,t'_5\ncr_H=H/ \ncl t_1,t_2,t_3,t_4,t_5,t'_5 \ncr_G=1$.
Therefore, $\tilde M_{56}$ is a complement in the 4-sphere.
\end{example}

\begin{example}
\label{m1091}
The side-pairing of manifold $M_{1091}$ is encoded by 53RR35.
All the generators are orientation reversing.
The presentation of $G$ is given by the following relations, which are grouped by row of diagrams like
in \figref{ridgecycles}.
\begin{displaymath}
\begin{array}{rrr}
\left\{
\begin{array}{r}
g^{-1}jh^{-1}i=1\\
g^{-1}dh^{-1}c=1\\
e^{-1}jf^{-1}i=1\\
e^{-1}df^{-1}c=1\\
hjgi=1\\
hc^{-1}gd^{-1}=1\\
f^{-1}je^{-1}i=1\\
f^{-1}c^{-1}e^{-1}d^{-1}=1
\end{array}
\right.
&
\left\{
\begin{array}{r}
g^{-1}lh^{-1}k=1\\
g^{-1}bh^{-1}a=1\\
e^{-1}lfk=1\\
e^{-1}bfa=1\\
hlgk=1\\
ha^{-1}gb^{-1}=1\\
fle^{-1}k=1\\
fa^{-1}e^{-1}b^{-1}=1
\end{array}
\right.
&
\left\{
\begin{array}{r}
i^{-1}k^{-1}ik=1\\
i^{-1}bia=1\\
c^{-1}k^{-1}d^{-1}k=1\\
c^{-1}bd^{-1}a=1\\
jlj^{-1}l^{-1}=1\\
ja^{-1}j^{-1}b^{-1}=1\\
dlcl^{-1}=1\\
da^{-1}cb^{-1}=1
\end{array}
\right.
\end{array}
\end{displaymath}
The table below lists the translations of interest to us in the boundary components.
The boundary components of $\tilde M_{1091}$ are AAAAAAAAA, making it a complement of 9 tori.
The boundary component $E_6$ is the one corresponding to vertex $v_{+++-}$, for which
the search for translations is the same as for $v_{++++}$.
\noindent
\begin{center}
\begin{tabular}{c|c|c|c|c}
comp. & type & translations in $8Q_C$ or $64 Q_C$ & holonomy &  self-conj?\\
\hline\rule[15pt]{0pt}{0pt}%
$E_1$ & $A$ & $c^{-1}h$, $a^{-1}h$, $c^{-1}b$ & trivial & no\\
$E_2$ & $H$ & $e^{-1}j$ & $x_3$ & no\\
$E_3$ & $H$ & $e^{-1}l$ &  $x_2$ & no\\
$E_4$ & $G$ & $i^{-1}k^{-1}$ &  $x_1$ & yes\\
$E_5$ & $A$ & $a^{-1}k, c^{-1}i$, $e^{-1}g$ & trivial & no\\
$E_6$ & $A$ & $a^{-1}l$, $c^{-1}j$, $e^{-1}h$ & trivial & no
\end{tabular}
\end{center}
Choose the translations:
\begin{displaymath}
\begin{array}{llllll}
t_1=c^{-1}b, &  t_2=e^{-1}j,  & t_3=e^{-1}l,  & t_4=i^{-1}k^{-1},   & t_5=c^{-1}i, & t_6=e^{-1}h,\\
t'_1=ct_1c^{-1},  & & & & t'_5=ct_5c^{-1},  & t'_6=ct_6c^{-1}.
\end{array}
\end{displaymath}
Due to self-conjugacy of $t_4$, $\ncl t_1,t_4,t_5,t_6\ncr_G=\ncl t_1,t'_1,t_4,t_5,t'_5,t_6,t'_6\ncr_H$.
In what follows, we refer to the relations above as entries in an $8\times 3$ matrix.
In the group $G/\ncl t_1,t_4,t_5,t_6\ncr $ we have $k^{-1}=i=c=b$ and $h=e$, so equation 43
gives $d=a$ and equation 23 gives $b=a^{-1}$.
Equations 32 and 42 convert to $l=ea^{-1}f^{-1}$ and $a^{-1}=ea^{-1}f^{-1}$, so $l=a^{-1}$.
Similarly, comparing 11 and 21 gives $j=a$.
The presentation of $G/\ncl t_1,t_4,t_5,t_6\ncr $ then reduces to generators $a$,$e$,$g$ and $f$ and
the relations, grouped as above (third column reduces to trivial relations):
\begin{displaymath}
\begin{array}{rr}
\left\{
\begin{array}{r}
g^{-1}ae^{-1}a^{-1}=1\\
e^{-1}af^{-1}a^{-1}=1\\
eaga^{-1}=1\\
f^{-1}ae^{-1}a^{-1}=1
\end{array}
\right.
&
\left\{
\begin{array}{r}
g^{-1}a^{-1}e^{-1}a=1\\
e^{-1}a^{-1}fa=1\\
ea^{-1}ga=1\\
fa^{-1}e^{-1}a=1
\end{array}
\right.
\end{array}
\end{displaymath}
In this array equations 11 and 41 give $f=g$, equations 22 and 32 imply $e=e^{-1}$,
and equations 12 and 42 give $g=f^{-1}$, further simplifying the presentation to
\begin{displaymath}
\left< a,e,f\ |\ e^2=f^2=1,\  aea^{-1}=f,\  afa^{-1}=e \right>\cong (\z_2 * \z_2)\rtimes \z,
\end{displaymath}
where $\z$ acts on $\z_2*\z_2$ by swapping the generators.
One can now show (next paragraph) that
$H/\ncl t_1,t_4,t_5,t_6\ncr\cong \z\oplus \z$, where the generators are $[t_2]=e^{-1}a$ and
$[t_3]=e^{-1}a^{-1}$.
Let $H_1=H/\ncl t_1,t'_1,t_4,t_5,t'_5,t_6,t'_6\ncr_H$.
We now have
\begin{multline*}
H/\ncl t_1,t'_1,t_2,t_3,t_4,t_5,t'_5,t_6,t'_6\ncr_H=\\
\left(H/\ncl t_1,t'_1,t_4,t_5,t'_5,t_6,t'_6\ncr_H\right)/\ncl [t_2],[t_3] \ncr_{H_1}=\\
\left(H/\ncl t_1,t_4,t_5,t_6\ncr \right)/\ncl [t_2],[t_3] \ncr_{H_1}=\\
\left( \left<[t_2]\right> \oplus \left< [t_3]\right>\right) /\ncl [t_2],[t_3] \ncr_{H_1}=1.
\end{multline*}
To show that $H/\ncl t_1,t_4,t_5,t_6\ncr\cong \z\oplus \z$ one can use either the explicit isomorphism
of $G/\ncl t_1,t_4,t_5,t_6\ncr$ with $(\z_2 * \z_2)\rtimes \z$, or apply the Reidemeister-Schreier method
to the presentation above using $\{1,e\}$ as the transversal.
In the latter case, we get the presentation 
\begin{displaymath}
\left< a_0,f_0,a_1,f_1\ | \ 
f_0f_1=a_0a_1^{-1}f_0^{-1}=a_0f_1a_1^{-1}=a_1a_0^{-1}f_1^{-1}=a_1f_0a_0^{-1}=1\right>,
\end{displaymath}
where $a_0=ae^{-1}$, $f_0=fe^{-1}$,
$a_1=ea$, and $f_1=ef$.
This easily reduces to $\left< a_0, a_1\, |\, a_0a_1=a_1a_0 \right>$.
We also note $[t_2]=e^{-1}a=ea=a_1$ and $[t_3]=e^{-1}a^{-1}=(ae^{-1})^{-1}=a_0^{-1}$.
\end{example}

\begin{example}
\label{m36}
The side-pairing of manifold $M_{36}$ is encoded by 1468AF.
The orientation-reversing generators are $e$,$f$,$i$,$j$,$k$,$l$.
The presentation of $G$ is given by the following relations.
\begin{displaymath}
\begin{array}{rrr}
\left\{
\begin{array}{r}
g^{-1}j^{-1}g^{-1}i=1\\
g^{-1}c^{-1}gc=1\\
e^{-1}jf^{-1}i=1\\
e^{-1}cfc=1\\
h^{-1}j^{-1}h^{-1}i=1\\
h^{-1}d^{-1}hd=1\\
ej^{-1}fi^{-1}=1\\
e^{-1}dfd=1
\end{array}
\right.
&
\left\{
\begin{array}{r}
g^{-1}l^{-1}h^{-1}k=1\\
g^{-1}a^{-1}ha=1\\
e^{-1}le^{-1}k=1\\
e^{-1}b^{-1}ea=1\\
g^{-1}kh^{-1}l^{-1}=1\\
gb^{-1}h^{-1}b=1\\
f^{-1}k^{-1}f^{-1}l^{-1}=1\\
f^{-1}b^{-1}fa=1
\end{array}
\right.
&
\left\{
\begin{array}{r}
i^{-1}l^{-1}i^{-1}k=1\\
i^{-1}b^{-1}ia=1\\
c^{-1}ld^{-1}k=1\\
c^{-1}a^{-1}da=1\\
jkjl^{-1}=1\\
ja^{-1}j^{-1}b=1\\
c^{-1}kd^{-1}l=1\\
cb^{-1}d^{-1}b=1
\end{array}
\right.
\end{array}
\end{displaymath}
\noindent
The following table identifies boundary component types and translations of interest.
The boundary components of $\tilde M_{36}$ are AAAABBAA making it a complement
of 6 tori and 2 Klein bottles.
\begin{center}
\begin{tabular}{c|c|c|c|c}
comp. & type & translations in $8Q_C$ or $64 Q_C$ & holonomy & self-conj?\\
\hline\rule[15pt]{0pt}{0pt}%
$E_1$ & $A$ & $c^{-1}$, $g^{-1}$ & trivial & no\\
$E_2$ & $A$ & $a^{-1}$, $e^{-1}j$ & trivial & no\\
$E_3$ & $J$ & none & $x_1$, $x_2$, $x_1x_2$ & n/a\\
$E_4$ & $J$ & none & $x_1$, $x_2$, $x_1x_2$ & n/a\\
$E_5$ & $A$ & $c^{-1}i^{-1}eg$, $a^{-1}k^{-1}ci$, $a^{-1}k^{-1}eg$ & trivial & no
\end{tabular}
\end{center}
Choose the translations:
\begin{displaymath}
\begin{array}{lllll}
t_1=a^{-1}b, & t_2=e^{-1}j, &  t_3= e^{-1}f^{-1},   & t_4=i^{-1}j^{-1},  & t_5=a^{-1}k^{-1}ci, \\
t'_1=kck^{-1},       & t'_2=eae^{-1},   & & & t'_5=k(c^{-1}i^{-1}eg)k^{-1}.
\end{array}
\end{displaymath}
Note that for boundary component $E_1$, although there are two translations in $8Q_C$, we choose
$a^{-1}b\in\ker\phi$ because neither $c$ nor $g$ is self-conjugate.
Note also that in boundary components $E_3$ and $E_4$ there are no translations in $8Q_C$ so
we are forced to choose translations in $\ker \phi$.
Lifts of $E_3$ and $E_4$ are of type $B$, so care must be taken that the translations are primitive
in the lifts.
In other words, the translations must be in the planes of rotation (in this case, the plane $x_1x_2$).

Since $t_3$ and $t_4$ are normal in nonorientable $E_3$ and $E_4$, they are self-conjugate,
so $\ncl t_3, t_4\ncr_H=\ncl t_3, t_4\ncr_G$.
In $G/ \ncl t_3,t_4\ncr$ we have $f=e^{-1}$ and $j=i^{-1}$.
Equations 42 and 82 then give $eae^{-1}=b$, $ebe^{-1}=a$, so $e(a^{-1}b)e^{-1}=b^{-1}a=(a^{-1}b)^{-1}$,
that is, $et_1e^{-1}=t_1^{-1}$.
Equation 31 reads as $e^{-1}i^{-1}ei=1$, so $et_2e^{-1}=t_2$ because $e$ and $j=i^{-1}$ commute.
Then $ \ncl t_1,t_2,t_3,t_4\ncr_H=\ncl t_1,t_2,t_3,t_4\ncr_G$.

In $G/ \ncl t_1,t_2,t_3,t_4\ncr$ we have $a=b$, which gives $eae^{-1}=a$, implying
\begin{displaymath}
\ncl t_1,t_2,t'_2,t_3,t_4\ncr_H=\ncl t_1,t_2,t'_2,t_3,t_4\ncr_G.
\end{displaymath}
In $G/\ncl t_1,t_2,t'_2,t_3,t_4\ncr$ we then have $a=1$ which implies $b=1$, $h=g$ and $d=c$.
Equation 32 gives $l=ek^{-1}e$.
Simplifying the presentation, we get
\begin{multline*}
G/\ncl t_1,t_2,t'_2,t_3,t_4\ncr =\left< c,e,g,k\ | \ egeg=g^{-1}c^{-1}gc=\right.\\
\left. e^{-1}ce^{-1}c=kgekeg=k^2=c^{-1}ekec^{-1}k=1\right>.
\end{multline*}
Unfortunately, the remaining translations $t'_1=kck^{-1}$, $t_5=kce^{-1}$ and $t'_5=kc^{-1}e^2gk^{-1}$
cannot easily be shown to be self-conjugate in this group, so we resort to the
Reidemeister-Schreier method to get $H/ \ncl t_1,t_2,t'_2,t_3,t_4\ncr$ before we continue.
Using $\{1,k\}$ as the transversal, it is not hard to see that the presentation for
$H/ \ncl t_1,t_2,t'_2,t_3,t_4\ncr$ has generators $c_0=c$, $g_0=g$, $e_0=ek^{-1}$,
$c_1=kck^{-1}$, $g_1=kgk^{-1}$, $e_1=ke$ and relations
\begin{displaymath}
\begin{array}{rrrrr}
e_0g_1e_1g_0=1 & g_0^{-1}c_0^{-1}g_0c_0=1 & e_1^{-1}c_1e_0^{-1}c_0=1 &
g_1e_1^2g_0=1 &  c_0^{-1}e_0^2c_1^{-1}=1\\
e_1g_0e_0g_1=1 & g_1^{-1}c_1^{-1}g_1c_1=1 & e_0^{-1}c_0e_1^{-1}c_1=1 &
g_0e_0^2g_1=1 & c_1^{-1}e_1^2c_0^{-1}=1
\end{array}
\end{displaymath}
The translations take the form $t'_1=c_1$, $t_5=c_1e_0^{-1}$, $t'_5=c_1^{-1}e_1e_0g_1$.
Adding $t'_1=t_5=t'_5=1$ to the above relations immediately gives $c_1=e_0=1$, $g_1=e_1^{-1}$.
Relation~15 then reduces to $c_0=1$ so~13 now says $e_1=1$, hence $g_1=1$.
Then $g_0=1$ follows from relation~21, therefore, $H/\ncl t_1,t'_1,t_2,t'_2,t_3,t_4,t_5,t'_5\ncr_H=1$.

In this example,  infinitely many other choices for $t_5$ and $t'_5$ also give the 4-sphere.
If we only add $c_1=1$ to the above equations, we get $c_0=e_0^2$, $e_1=e_0$ and
$g_1=e_0^{-2}g_0^{-1}$, and the presentation simplifies to
$\left<e_0, g_0\,|\,e_0g_0e_0^{-1}g_0^{-1}=1\right>$.
Thus, $H/\ncl t_1,t'_1,t_2,t'_2,t_3,t_4\ncr_H=\left<e_0\right>\oplus \left<g_0\right>$.
The basis elements of $\pi_1 \tilde E_5=\z^3$ are $v_1=c^{-1}i^{-1}eg$, $v_2=a^{-1}k^{-1}ci$,
and $v_3=a^{-1}k^{-1}eg$; they are sent to $g_0$, $e_0^{-1}$ and $e_0g_0$.
The conjugates $v'_i=kv_i k^{-1}$ generate $\pi_1\tilde E'_5$; they are sent
to $g_0^{-1}$, $e_0$ and $e_0^{-1}g_0^{-1}$.
Any primitive translations in $\tilde E_5$, $\tilde E'_5$ then have form  $t_5=pv_1-qv_2+rv_3$ or
$t'_5=-p' v'_1+q' v'_2-r' v'_3$ where $\gcd(p,q,r)=\gcd(p',q',r')=1$.
If $\psi:H\to H/\ncl t_1,t'_1,t_2,t'_2,t_3,t_4\ncr_H$ is the quotient map, then $\psi(t_5)=(q+r, p+r)$
and $\psi(t'_5)=(q'+r',p'+r')$, as written in terms of the basis $\{e_0, g_0\}$.
Then 
\begin{displaymath}
H/\ncl t_1,t'_1,t_2,t'_2,t_3,t_4,t_5,t'_5\ncr_H=1
\Longleftrightarrow
\left|
\begin{array}{cc}
q+r & q'+r'\\
p+r & p'+r'
\end{array}
\right|=1.
\end{displaymath}
Infinitely many choices of $p$, $q$, $r$, $p'$, $q'$, $r'$ exist that satisfy the last condition.
For example, settting $r=r'=0$, it is easy to see that the resulting equation
$qp'-pq'=1$ can be satisfied for an arbitrary relatively prime pair of numbers $p$ and $q$
once suitable $q'$ and $r'$ are found ($q'$ and $r'$ are consequently relatively prime).

\end{example}

\begin{example}
\label{otherex}
For the remaining examples from Theorem~\ref{mainthm}, the proofs
proceed in a similar way.  In \tableref{table2} we list the
choices for translations that make the orientable double cover a
complement in $S^4$.  Types of boundary components are identified in
the order we have used here, note that this order is different from
tables in \cite{Ratcliffe-Tschantz}.  The translations are listed in
the order one sets them equal to 1.  The letters RS, if applicable,
stand in the place where we were forced to pass to a presentation of
$H/\ncl \text{translations}\ncr$ via the Reidemeister-Schreier method.
Examples that do not employ it are typically straightforward.
Appearances of translations labeled $t'$ indicate that the
corresponding boundary component lifted to two components in the
orientable double cover.\end{example}

\begin{table}[ht!]\small\anchor{table2}
\begin{center}
\begin{tabular}{c|c|c|c}
number & side-pairing & boundary & translation choice\\
\hline\rule[15pt]{0pt}{0pt}%
1011 & 14FF28 & GGGGG & $t_1=c$, $t_2=a$, $t_3=k$,  $t_4=i$, \\
&&& $t_5=e^{-1}g$ (see \cite{Ivansic3})\\
\hline\rule[15pt]{0pt}{0pt}%
71 & 13EB34 & GGHGA & $t_1=a^{-1}h$, $t_2=a$, $t_4=k$, \\
&&& $t_5=c^{-1}i$, $t'_5=ct_5c^{-1}$, $t_3=e^{-1}l$ \\
\hline\rule[15pt]{0pt}{0pt}%
23 & 1569A4 & GAGAH & $t_1=c^{-1}h$, $t_2=a$, $t'_2=ct_2c^{-1}$,\\
&&& $t_3=e^{-1}d^{-1}$, $t_4=k$,\\
&&&  $t'_4=ct_4c^{-1}$, $t_5=c^{-1}i^{-1}eg$  \\
\hline\rule[15pt]{0pt}{0pt}%
1092 & 53FFCA & AGGAGG & $t_1=c^{-1}h$, $t'_1=ct_1c^{-1}$, $t_2=a^{-1}i^{-1}$, \\
&&& $t_3=c^{-1}k^{-1}$, $t_4=i^{-1}h^{-1}$, $t'_4=ct_4c^{-1}$ \\
&&& $t_5=e^{-1}g$, $t_6=e^{-1}h$\\
\hline\rule[15pt]{0pt}{0pt}%
231 & 1569F4 & GBGBH & $t_1=c^{-1}h$, $t_2=a$, $t'_2=ct_2c^{-1}$,\\
&&& $t_3=e^{-1}d^{-1}$, $t_4=k$, $t'_4=ct_4c^{-1}$,\\
&&&  $t_5=a^{-1}k^{-1}ci^{-1}eg$  \\
\hline\rule[15pt]{0pt}{0pt}%
112 & 13C874 & GGHBG & $t_1=g$, $t_2=a$, $t_4=k$, $t'_4=ct_4c^{-1}$, \\
&&& $t_5=c^{-1}ieg$, $t_3=e^{-1}d^{-1}l$\\
\hline\rule[15pt]{0pt}{0pt}%
92 & 1348EC & GAHIH & $t_1=g$, $t_4=k^{-1}l^{-1}$,\\
&&& $t_2=e$, $t_3=e^{-1}d^{-1}l$,\\
&&& RS, $t'_2=cac^{-1}$, $t_5=c^{-1}jak$\\
\hline\rule[15pt]{0pt}{0pt}%
51 & 156A9C & HGABG & $t_2=a$, $t_3=c^{-1}l$, $t'_3=ct_3c^{-1}$,\\
&&& $t_4=i^{-1}h^{-1}$, $t'_4=ct_4c^{-1}$,\\
&&& $t_5=c^{-1}i^{-1}f^{-1}g$, $t_1=c^{-1}ag^{-1}$\\
\hline\rule[15pt]{0pt}{0pt}%
40 & 143CF9 & GAGAJ & $t_1=c$, $t_3=e^{-1}k^{-1}$, $t_5=e^{-1}h^{-1}fg$,\\
&&& $t_4=g^{-1}l$, $t_2=a$,\\
&&& RS, $t'_2=i(e^{-1}j)i^{-1}$, $t'_4=i(i^{-1}h)i^{-1}$\\
\end{tabular}
\end{center}
\nocolon\caption{}\label{table2}
\end{table}

This completes the proof of Theorem~\ref{mainthm}.
\end{proof}

\eject

We conclude with some related examples.

\begin{example}
\label{m1091supp}
(Examples of complements in simply-connected manifolds with higher Euler characteristic)\qua
It was shown in \cite{Ivansic3} that index-$n$ cyclic covers of $\tilde M_{1011}$ are complements
of $4n+1$ tori in a simply-connected closed manifold $N$ with Euler characteristic $2n$.
We show that many other manifolds from Theorem~\ref{mainthm} have an analogous property.
Two propositions from~\cite{Ivansic3} are used for this purpose.

Let $M$ be a hyperbolic manifold with boundary components $E_1,\dots,E_m$ and let
$H=\pi_1 M$, $K_i=\pi_1 E_i$.
If we take the cover corresponding to a normal subgroup $H_0$ of $H$,
\mbox{\cite{Ivansic3} Proposition~2.3} asserts that the number of path-components of $p^{-1}(E_i)$
is equal to the index of the image of $K_i$ in $H/H_0$ under the quotient map.
Now let $t_i\in \pi_1 E_i$ be a translation normal in $E_i$, $i=1,\dots,m$.
Setting $H_0=\ncl t_1,\dots,t_m\ncr_H$, suppose that $H/H_0$ is a finite group of order $l$, and that
$t_i$ is primitive $K_i \cap H_0$.
Then \mbox{\cite{Ivansic3} Proposition~2.4} states that the index-$l$ cover of $M$ corresponding to
subgroup $H_0$ is a complement inside a simply-connected closed manifold.

From example~\ref{m56} it is easily seen that
$H/\ncl t_1,t_2,t'_2,t_3,t_4,t_5\ncr_H=\z=\left<gc\right>=\left<[t'_5]\right>$.
It follows that $H/\ncl t_1,t_2,t'_2,t_3,t_4,t_5, (t'_5)^m\ncr_H=\z_m$, so the above discussion
shows that an $m$-fold cover of  $\tilde M_{56}$ is a complement in a simply-connected 4-manifold
$N$ with $\chi(N)=2m$.
To find the link type, consider the generators of boundary components of $\tilde M_{56}$ and
their images under the quotient map $H\to H/\ncl t_1,t_2,t'_2,t_3,t_4,t_5\ncr_H$:
\begin{center}
\begin{tabular}{c|c|c|c}
component & type & generators  & images in $\z=\left<gc\right>$ \\
\hline\rule[15pt]{0pt}{0pt}%
$\tilde E_1$ & B &$a^{-1}b$, $g^{-1}h$, $c^{-1}h$ & $1$, $g^{-2}=(gc)^{-2}$, $gc$\\
$\tilde E_2$ & B & $a$, $e$, $i^{-1}j$ & $1$, $1$, $1$\\
$\tilde E'_2$ & B &$cac^{-1}$, $cec^{-1}$, $c(i^{-1}j)c^{-1}$& $1$, $1$, $1$ \\
$\tilde E_3$ & A & $c^{-2}$, $e^{-1}dl$, $e^{-1}f$ & $1$, $1$, $1$ \\
$\tilde E_4$ & A & $g^{-1}h^{-1}$, $i^{-2}$, $i^{-1}k^{-1}$ & $1$, $1$, $1$ \\
$\tilde E_5$ & A & $a^{-1}k^{-1}a^{-1}k$, $a^{-1}k^{-1}eg$, $c^{-1}i$ & $1$, $gc$, $1$\\
$\tilde E'_5$ & A & $c(a^{-1}k^{-1}a^{-1}k)c^{-1}$, & $1$, $gc$, $1$\\
&& $c(a^{-1}k^{-1}eg)c^{-1}$, $c(c^{-1}i)c^{-1}$ &\\ 
\end{tabular}
\end{center}
The results quoted above imply that the number of components in the lift of a boundary
component $\tilde E$ is simply the order of the group $\z/\left<A, [t'_5]^m\right>$, where $A$ is the
set of images of the generators of $\pi_1 \tilde E$.
Those orders are, respectively, 1, $m$, $m$, $m$, $m$, 1, 1.

Let $H_0=\ncl t_1,t_2,t'_2,t_3,t_4,t_5, (t'_5)^m\ncr_H$.
Note that $(t'_5)^m$ is primitive in $H_0\cap \pi_1\tilde E'_5$, since  powers smaller than $m$ of
$t'_5$ are not in $H_0$ because they map to nontrivial elements of $H/H_0=\z_m$.

The only generator of $\pi_1 \tilde E_1$ with nontrivial holonomy is $c^{-1}h$, mapping under
the restriction $\psi: \pi_1\tilde E_1\to \z_m$ to the generator of $\z_m$.
If $m$ is even, elements of $\ker\psi$ will all have trivial holonomy, meaning that $\tilde E_1$ lifts
to a collection of manifolds of type $A$; if $m$ is odd, some elements of $\ker\psi$ will have nontrivial
holonomy, so $\tilde E_1$ will lift to a collection of manifolds of type $B$.

Putting everything together, we see that an $m$-fold cyclic cover of $\tilde M_{56}$ is a complement 
inside a simply-connected closed manifold $N$ with $\chi(N)=2m$ of
\begin{center}
\begin{tabular}{l}
$2m+3$ tori and $2m$ Klein bottles, if $m$ is even;\\
$2m+2$ tori and $2m+1$ Klein bottles, if $m$ is odd.
\end{tabular}
\end{center}

At the end of Example~\ref{m36} we saw that $H/\ncl t_1,t'_1,t_2,t'_2,t_3,t_4\ncr_H=\z\oplus\z$,
where $\pi_1\tilde E_5$ and $\pi_1 \tilde E'_5$ map surjectively to $\z\oplus\z$ under the
quotient map.
This now allows for many choices of $t_5$ and $t'_5$ that make
$H/\ncl t_1,t'_1,t_2,t'_2,t_3,t_4,t_5,t'_5\ncr_H$ finite.
With notation from Example~\ref{m36} it is clear that 
\begin{displaymath}
H/\ncl t_1,t'_1,t_2,t'_2,t_3,t_4,t_5,t'_5\ncr_H\text{ is finite}
\Longleftrightarrow
\left|
\begin{array}{cc}
q+r & q'+r'\\
p+r & p'+r'
\end{array}
\right|\ne0.
\end{displaymath}
Note that $t_5$ and $t'_5$ do not have to be primitive in $\pi_1 \tilde E_5$ or $\pi_1 \tilde E'_5$,
so there is no condition on $\gcd(p,q,r)$ or $\gcd(p',q',r')$.
It turns out $t_5$ and $t'_5$ are always primitive in $H_0\cap \pi_1 \tilde E_5$ and 
$H_0\cap \pi_1 \tilde E'_5$, where $H_0=\ncl t_1,t'_1,t_2,t'_2,t_3,t_4,t_5,t'_5\ncr_H$.

The number of boundary components of the cover of $\tilde M_{36}$ corresponding to $H_0$ is
determined by considering the order of the groups
\begin{displaymath}
\z\oplus\z/\left< (q+r,p+r), (q'+r',p'+r'), A\right>,
\end{displaymath}
where $A$ is the set of images of generators
of boundary components of $\tilde M_{56}$.
Care must be taken to ascertain whether boundary components of type B lift to type B or A.

Among the many possibilities here, we just consider the case $q=p'=r=r'=0$, with $p$ and $q'$
arbitrary.
Under those conditions, one can see that $\tilde M_{36}$ has a finite cover with deck group
$\z_p\oplus \z_{q'}$ that is the complement in a simply-connected manifold $N$ with $\chi(N)=2pq'$~of
\begin{center}
\begin{tabular}{l}
$(1+2p)\gcd(2,q')+\gcd(2p,q')+2+2p$ tori, if $q'$ is even;\\
$(1+2p)\gcd(2,q')+\gcd(2p,q')+2$ tori and $2p$ Klein bottles, if  $q'$ is odd.
\end{tabular}
\end{center}

The table below lists other examples of complements in a simply-connected manifold $N$ that we
found using the same method.
The last column indicates which translations different from the ones in Example~\ref{otherex} are set
equal to~1 to cause $H/\ncl \text{translations}\ncr_H$ to be a finite group.
The cover is then the one corresponding to the kernel of $H\to H/\ncl \text{translations}\ncr_H$.
The Euler characteristic of $N$ is always twice the cardinality of the group of deck transformations.
\begin{center}
\begin{tabular}{c|c|l|c}
A cover  & with & is a complement in a simply- &  Altered choice\\
of mfd. & deck group & connected closed manifold of: & of translations\\
\hline
\rule[15pt]{0pt}{0pt}%
$\tilde M_{1091}$ & $\z_m\oplus\z_n $ & $4\gcd(m,n)+5$ tori & $t_2^m$ and $t_3^n$\\
\hline
\rule[15pt]{0pt}{0pt}%
$\tilde M_{71}$ & $\z_m$ & $3m+3$ tori  & $t_3^m$\\
\hline
\rule[15pt]{0pt}{0pt}%
$\tilde M_{112}$ & $\z_m$ & $3m+3$ tori, $m$ even & $t_3^m$\\
&& $3m+1$ tori, 2 K. bottles, $m$ odd &\\
\hline
\rule[15pt]{0pt}{0pt}%
$\tilde M_{92}$ & $\z_m$ & $4m+2$ tori, $m$ even & $t_5^m$ \\
&& $4m+1$ tori, 1 K. bottle, $m$ odd &\\
\hline
\rule[15pt]{0pt}{0pt}%
$\tilde M_{92}$ & $\z_m$ & $3m+3$ tori, $m$ even & $(t'_2)^m$ \\
&& $3m+2$ tori, 1 K. bottle, $m$ odd &\\
\hline
\rule[15pt]{0pt}{0pt}%
$\tilde M_{40}$ & $\z_m$ & $2m+4$ tori, $m$ K. bottles & $(t'_2)^m$ or $(t'_4)^m$\\
\end{tabular}
\end{center}
\end{example}

\begin{example}
\label{aspherical}
(More examples of aspherical homology spheres.)
Ratcliffe and Tschantz produced the first examples (infinitely many) of aspherical 4-manifolds that are
homology spheres.
The construction in~\cite{Ratcliffe-Tschantz2} used $\tilde M_{1011}$: by choosing fibers $t_i$ of
boundary components $E_i$
as in Example~\ref{otherex} and filling them in with a disc, one gets a manifold $N$ with
$\pi_1 N=1$, thus $H_1 N=0$.
By Proposition~\ref{crithom}, $H_1 \tilde M_{1011}=\z^5$, generated by the fibers
$t_1,\dots,t_5$.
If the fibers are altered,  $\pi_1 N$ will likely fail to be trivial.
However, one can retain $H_1 N=0$ if $t_i$ is replaced with $t_i+s_i$, where $s_i$ is a carefully
chosen element of $H_1 E_i=\z^3$ that can have arbitrarily large length.
When fibers are chosen to be sufficiently long, the manifold $N$ supports a metric of nonpositive
curvature and is therefore aspherical ($N$ is a homology sphere because $H_1 N=0$
and $\chi(N)=2$).
We note that this construction goes through for any of the orientable double covers
from~\ref{mainthm} whose boundary components are all 3-tori (i.e. which are
complements of a collection of tori in $S^4$), namely
$\tilde M_{71}$, $\tilde M_{23}$, $\tilde M_{1092}$ and $\tilde M_{1091}$.
\end{example}

\Addresses\recd

\end{document}